\documentclass{article}
\usepackage[utf8]{inputenc}
\usepackage{amsmath,amsthm,amssymb,amsfonts,color,bbm}
\usepackage{tikz}
\usetikzlibrary{calc}
\usepackage{diagbox}
\usepackage{authblk}
\usepackage{subcaption}
\usepackage{hyperref}
\newtheorem{theorem}{Theorem}

\newtheorem{lemma}{Lemma}
\newtheorem{corollary}{Corollary}
\newtheorem{remark}{Remark}

\newtheorem{assumption}{Assumption}
\newcommand{\mathrom}{\mathsf}
\newcommand{\softmax}{\textrm{Softmax}}
\title{Irrationality as a mean of regularization in Bayesian Persuasion}
\author{Romain Duboscq\footnote{romain.duboscq@math.univ-toulouse.fr}}
 \affil{Institut de Math\'ematiques de Toulouse ; UMR5219\\Universit\'e de Toulouse ; CNRS\\INSA, F-31077 Toulouse, France}
\author{Frederic de Gournay\footnote{degourna@insa-toulouse.fr}}
 \affil{Institut de Math\'ematiques de Toulouse ; UMR5219\\Universit\'e de Toulouse ; CNRS\\INSA, F-31077 Toulouse, France}
\begin{document}

\maketitle
\begin{abstract}
     We study a regularized variant of the Bayesian Persuasion problem, where the receiver's decision process includes a divergence-based penalty that accounts for deviations from perfect rationality. This modification smooths the underlying optimization landscape and mitigates key theoretical issues, such as measurability and ill-posedness, commonly encountered in the classical formulation. It also enables the use of scalable second-order optimization methods to compute numerically the optimal signaling scheme in a setting known to be NP-hard. We present theoretical results comparing the regularized and original models, including convergence guarantees and structural properties of optimal signaling schemes. Analytical examples and numerical simulations illustrate how this framework accommodates complex environments while remaining tractable and robust. A companion Python library, \textit{BASIL}\footnote{https://plmlab.math.cnrs.fr/degourna/bayesian-persuasion-by-regularization}, makes use of all the practical insights from this article.
\end{abstract}
\newpage
\tableofcontents

\section{Introduction}

The Bayesian Persuasion framework, introduced by Kamenica and Gentzkow \cite{kamenica2011bayesian}, provides a foundational model for understanding how informed agents can strategically design information structures to influence the actions of less informed receivers. By framing persuasion as the choice of a signal structure that shapes the receiver's posterior beliefs, the model has offered deep insights into information transmission in economics, political science, and beyond. We refer the reader to \cite{kamenica2019bayesian,bergemann2019information} and the references therein for the different models and applications.

Despite its elegance, the classical Bayesian Persuasion model assumes fully rational agents and frictionless belief updating, often limiting its applicability in environments where human or institutional behavior deviates from ideal Bayesian reasoning. In many practical settings, receivers exhibit bounded rationality, behavioral biases, or computational limitations that affect how they process and respond to information. This issue has been investigated in different contexts where the agents have different priors \cite{alonso2016bayesian,galperti2019persuasion}, make a non-Bayesian updating \cite{epstein2006axiomatic,epstein2008non,augenblick2021belief,de2022non,zhao2022pseudo}, are rationally inattentive \cite{bloedel2018persuasion,wei2021persuasion,lipnowski2020attention} or have a prior bias \cite{lee2023cheap}.

In this work, we extend the canonical Bayesian Persuasion model by introducing a regularization term (specifically, a divergence) in the receiver's optimization problem. This regularization induces a smoother, more tractable optimization landscape, providing both analytical clarity and computational robustness. To be more specific, this modification serves two main purposes: it models potential irrationality or sub-optimal behavior on the part of the receiver, and it addresses both theoretical and computational challenges inherent in the standard formulation. On the theoretical side, as we will show, the classical model may lead to measurability issues and difficulties in identifying the set of admissible actions. On the computational side, the persuasion problem is known to be NP-hard, and prior approaches have attempted to overcome this complexity through specialized linear programming techniques (see, e.g., \cite{dughmi2016algorithmic, dughmi2017aalgorithmic, dughmi2017balgorithmic, babichenko2017algorithmic, dughmi2019hardness, gradwohl2022algorithms}). In contrast, our formulation enables the use of quasi-Newton methods, offering a more efficient and flexible computational framework.

In Section \ref{sec:problem}, we develop the mathematical framework for both the standard Bayesian Persuasion model and its regularized extension. Section \ref{sec:results} presents our main theoretical findings. For the classical model, we establish a lower bound on the number of messages required in an optimal signaling scheme. While this issue is typically addressed via the revelation principle, that approach does not directly apply when receivers possess private types. Regarding the regularized formulation, we demonstrate that it generally constitutes a sub-optimal approximation of the original problem. Nevertheless, under certain assumptions, we prove that it serves as a valid approximation in the sense that any sequence of minimizers converges, up to a subsequence, to a minimizer of the non-regularized problem. We also provide numerical insights and discuss the use of second-order optimization methods, which yield efficient and fast algorithms. Finally, in Section \ref{sec:Examples}, we illustrate the framework and support our theoretical arguments through both analytical examples, which admit closed-form solutions, and numerical simulations, which highlight the model's flexibility and applicability to more complex environments.

\section{Setting of the problem}\label{sec:problem}

\subsection{Notations and preliminaries}
For any $n\in\mathbb{N}$, we denote $\Delta_n$ the simplex of dimension $n$, that is
\begin{equation*}
    \Delta_n = \left\{x\in(\mathbb{R}^+)^n\text{ such that } \sum_{k = 1}^nx_k = 1\right\}.
\end{equation*}
For any measurable space $(\mathcal{E},\mathsf{E})$, where $\mathsf{E}$ is a $\sigma$-algebra of $\mathcal{E}$, we denote $\mathsf{P}(\mathcal{E})$ the set of probability measures on $(\mathcal{E},\mathsf{E})$. When $\mathcal{E}$ is finite, we will always endow it with the $\sigma$-algebra $\mathsf{E}=\mathrm{P}(\mathcal{E})$, the set of all parts of $\mathcal{E}$. Moreover, $\mathsf{P}(\mathcal{E})$ is a compact space which identifies as $\Delta_{|\mathcal{E}|}$ where $|\mathcal{E}|$ denotes the cardinal number of $\mathcal{E}$. When $(\mathcal{E},d)$ is a metric space, the space $(\mathcal{P}(\mathcal{E}),d_{\mathrm{Pr}})$ is a metric space \cite{billingsley2013convergence} where $d_{\mathrm{Pr}}$ is the Prokhorov metric given by
\begin{equation*}
    d_{\mathrm{Pr}}(\mu,\nu) = \inf\{r>0:\; \mu(E)\leq \nu(E_r) + r\;\text{and}\; \nu\leq \mu(E_r) + r,\;\forall E\in\mathsf{B}(\mathcal{E})\},
\end{equation*}
where $\emptyset_r = \emptyset$ as well as, for any $E\neq \emptyset$,
\begin{equation*}
    E_r = \{e\in\mathcal{E}:\; d(e,E)<r\}.
\end{equation*}
and $\mathsf{B}(\mathcal{E})$ is the Borel $\sigma$-algebra of $\mathcal{E}$.
Furthermore, if $(\mathcal{E},d)$ is compact, then the space $(\mathcal{P}(\mathcal{E}),d_{\mathrm{Pr}})$ is also compact.

\subsection{Agents and variables}
The agents are a sender (referred to as "she") and one or more receivers (referred to as "them", obnoxious to the effective number of receivers). We first describe the variables at play, for the sake of simplicity, we suppose that each considered set is finite.
\begin{itemize}
    \item The set of receivers is denoted $\mathcal L$ and is of cardinality $|\mathcal L| = L \ge 1$. The generic notation for a receiver is $\ell \in \mathcal L$ and it is common to identify  $\mathcal L$ with $[1,\dots,L]$.
    \item The set of states is denoted by $\mathcal S$. The states are denoted $s$ and there are $|\mathcal S| = S$ states in total. The states model a source of uncertainty that is common to all agents. 
    \item Each receiver $\ell$ has a type $t_\ell \in \mathcal T_l$.
    There are $|\mathcal T_{\ell}| = T_\ell \ge 1$ types available to the receiver $\ell$. The type of the receiver changes its utility and prior on $\mathcal S$. For each $t_\ell\in \mathcal T_\ell$, the probability of the $\ell$-th receiver to be of type $t_\ell$ is given by $\eta_\ell(t_\ell)$ with $\eta_\ell \in\mathrom{P}(\mathcal T_l)$. The type of a receiver is unknown to every other agents, however every agent is aware of $\eta_\ell$, the probability of being of a certain type. The set of types is denoted by $\mathcal T=\bigotimes_{\ell \in \mathcal L} \mathcal T_\ell$ and is of cardinal $T=\prod_{\ell \in \mathcal L}  T_\ell$.
    \item The set of messages is denoted by $\mathcal M$. These will play an important role since they are at the heart of the way the information from the sender is transmitted to the receiver (see below the communication policy, or signal). The generic notation for a message is $m\in \mathcal M$, there are $|\mathcal M| = M$ messages in total.
    \item Each receiver $\ell$ has to pick an action $a_\ell$ in his available set of actions $\mathcal A_\ell$ of cardinal $|\mathcal A_{\ell}| = A_\ell$. The chosen action is known to the other agents, hence the set of available actions do not depend on the type (or else the receiver would disclose his type). The set of actions is denoted $\mathcal{A}=\bigotimes_{\ell \in \mathcal L} \mathcal{A}_\ell$ and is of cardinal $A=\prod_{\ell \in \mathcal L} A_\ell$.
    \item Finally, each receiver $\ell$ or type $t_\ell$ has a utility $(s,a_\ell)\rightarrow u_{t_\ell}(s,a_\ell)$ that depends on the state $s\in S$ and the chosen action $a_\ell \in \mathcal A_\ell$. This utility drives the choice of action picked by the receiver. The sender has a utility $(s,a)\rightarrow v(s,a)$ that depends on the state $s\in S$ and on each chosen action $a=(a_\ell)_{\ell \in \mathcal L}\in \mathcal A$
\end{itemize}
\subsection{The receivers' problem : choosing the action}\label{sec:receiverproblem}
In this section, we focus on describing the process by which the receiver $\ell$ of type $t_\ell$ chooses his action.
For any prior $\nu \in \mathrom P(\mathcal S)$ on the states, the receiver of type $t_\ell$ maximises its utility and computes
\begin{equation}\label{eq:probinitreceiver:argmax}
    \Theta_{t_{\ell}}(\nu)=\underset{\theta\in \mathrom{P}(\mathcal A_\ell)}{\mathrm{argmax}}\;\sum_{(s,a_\ell)\in\mathcal S\times \mathcal{A}_\ell} u_{t_\ell}(s,a_\ell)\theta(a_\ell)\nu(s).
\end{equation}
The elements in $\Theta_{t_{\ell}}(\nu)$ are the "acceptable strategies for the receiver of type $t_{\ell}$ under $\nu$" and are {\em a priori } not unique. It is however well known that $\Theta_{t_{\ell}}(\nu)$ is the convex hull of dirac masses located on $\mathcal{A}^\star_{t_\ell}(\nu)$ the "admissible actions" (or "pure strategies") defined as
\begin{equation}\label{eq:probinitreceiver:admissibleactions}
    \mathcal{A}^\star_{t_\ell}(\nu)=\underset{a \in \mathcal A_\ell}{\mathrm{argmax}}\;\sum_{s\in\mathcal S} u_{t_\ell}(s,a)\nu(s).
\end{equation}
If there is only one admissible action, that is the cardinal of $\mathcal{A}^\star_{t_\ell}$ is one, then the receiver picks up this action. If there are several admissible actions, the receiver picks the actions that are the most favorable to the sender. As soon as they compute their admissible strategies, they make it public to the rest of the agents and hence each agent is aware of $\Theta_{t}(\nu)$ given by
$$
\Theta_{t}(\nu) = \left\{ a\in\mathcal{A}\mapsto \prod_{\ell\in\mathcal{L}}\theta_{t_\ell,\nu}(a_{\ell})\text{ such that }\theta_{t_\ell,\nu}\in \Theta_{t_{\ell}}(\nu)\;,\forall \ell\in\mathcal{L}  \right\}.
$$
Finally they choose a global strategy $\theta_{t,\nu}^\star\in \Theta_{t}(\nu)$ a solution to
\begin{equation}
\label{eq:compliant}
   \max_{\theta \in \Theta_{t}(\nu) }\left(\sum_{(s,a)\in\mathcal{S}\times\mathcal{A}} v(s,a)\nu(s)\theta(a)\right) ,
\end{equation}
where, $v$ is the utility of the sender and $\nu$ is her prior. Solutions to the linear programming problem \eqref{eq:compliant} may not be unique, but one can always select a solution that corresponds to a vertex of the feasible polytope. Consequently, there exists at least one choice of $\theta_{t,\nu}^\star$ that is a Dirac measure concentrated on a single action $a_t^\star(\nu)$, such that for every type $t_\ell$, we have $a_{t_\ell}^\star(\nu) \in \mathcal{A}_{t_\ell}^\star(\nu)$, i.e., an admissible action. By convention, we assume that receivers of type $t_\ell$ adopt such a (pure) strategy.

In other words, for every $t \in \mathcal{T}$ and $\nu \in \mathrom{P}(\mathcal{S})$, we define $a_t^\star(\nu) \in \mathcal{A}_t^\star(\nu)$, where $\displaystyle\mathcal A_t^\star(\nu)=\bigotimes_{\ell\in \mathcal L} \mathcal A_{t_\ell}^\star(\nu)$, as one solution to the following problem: 
\[\sum_{s \in \mathcal S} v(s,a^\star_t(\nu))\nu(s) \ge \sum_{s \in \mathcal S} v(s,a)\nu(s), \quad \forall a \in \mathcal A_t^\star(\nu).\]

\begin{remark}
Throughout this section, we have implicitly assumed that all agents share the same prior $\nu \in \mathcal{P}(\mathcal{S})$. However, one can also handle the case of heterogeneous priors using the following trick: select any prior $\nu \in \mathcal{P}(\mathcal{S})$ such that each agent's prior is absolutely continuous with respect to $\nu$, and, for each agent, scale its utility by the ratio of its own prior to $\nu$. The agent will behave the same way under this new utility and the shared prior $\nu$.
\end{remark}

\subsection{The sender's problem : choosing the message}\label{sec:senderproblem}
The sender is the main agent and she'll want to maximize her utility, denoted $v:\mathcal{S}\times \mathcal{A}_1\times\ldots\times\mathcal{A}_{L}\mapsto \mathbb{R}$, which depends on the state of the world and the actions from the receivers. We suppose that each agent share the same prior $\mu\in\mathrom{P}(\mathcal{S})$ which is enhanced by a message conveyed by the sender. Indeed, she is able to design a communication policy (or signal) $\pi = (\pi(m|s))_{(m,s)\in\mathcal{M}\times\mathcal{S}}$ in a way that will influence the actions of receivers to her benefit. The value $\pi(m|s)$ is to be understood as the probability of receiving the message $m$ given $s$, the state of the world. We observe that, for each state $s$, $\sum_{m\in \mathcal M} \pi(m|s)=1$ and we denote $\mathrom{P}_{\mathcal{S}}(\mathcal{M})$ the set of available communication policies 
\[\mathrom{P}_{\mathcal{S}}(\mathcal{M}) =\left\{\pi(m|s) \text{ such that } \pi(\cdot|s)\in \mathrom{P}(\mathcal{M}) \text{ for all } s\right\},\]
With a communication policy $\pi \in\mathrom{P}_{\mathcal{S}}(\mathcal{M})$ at hand, after receiving a message $m\in\mathcal{M}$, the receivers update their prior $\mu$ by Bayes' rule and compute their posterior
\begin{equation}
\label{eq:define:nu}
    \nu_{m,\pi}(s) = \frac{\pi(m|s)\mu(s)}{p(m)} \quad\text{ with }\quad p(m) = \sum_{\tilde{s}\in\mathcal{S}} \pi(m|\tilde{s})\mu(\tilde{s}),
\end{equation}
where $p(m)$ is the probability of receiving the message $m$. With this posterior $\nu_{m,\pi}$, the receivers then compute their acceptable strategies $\Theta_{t,\nu_{m,\pi}}$ and reveal them to the sender. The sender then chooses amongst the available strategies the most favorable ones.

From here, the sender needs to design her communication policy $\pi$ in order to maximize her utility. Let us notice that, when designing $\pi$, there are several random variables whose realizations she's not aware of: the state of the world, the message and the type of each receiver. The information on these variables is encoded through a distribution $\tilde \eta\in\mathrom{P}(\mathcal{S}\times\mathcal{M}\times \mathcal{T})$. Since the types are independent of the states and of the messages, we can factorize, for any $s\in\mathcal{S}$, $m\in\mathcal{M}$ and $t\in \mathcal{T}$,
\begin{equation*}
    \tilde{\eta}(s,m,t) = \tilde{\mu}(s,m) \eta(t) \text{ and } \eta(t)=\prod_{\ell \in \mathcal L}\eta_\ell(t_\ell),
\end{equation*}
where $\tilde{\mu}$ is computed thanks to her communication policy as well as her prior $\mu$. That is, we have, for any $(s,m)\in\mathcal{S}\times\mathcal{M}$,
\begin{equation*}
    \tilde{\mu}(s,m) = \pi(m|s)\mu(s).
\end{equation*}
In the end, when designing $\pi$ and by using \eqref{eq:define:nu}, the sender faces the following maximization problem
\begin{equation}
\label{eq:middleformulation}
    \max_{\pi\in\mathrom{P}_{\mathcal{S}}(\mathcal{M})}\left(\sum_{(m,t)\in\mathcal{M}\times\mathcal{T}}\max_{\theta \in \Theta_t(\nu_{m,\pi}) }\left(\sum_{(s,a)\in\mathcal{S}\times\mathcal{A}} v(s,a)\nu_{m,\pi}(s) \theta(a)\right)\eta(t)p(m) \right).
\end{equation}
For a given $t$ and $m$, the problem of maximization in $\theta$ has been discussed in Section~\ref{sec:receiverproblem} and the optimal solution is a Dirac located at the action $a_t^\star(\nu_{m,\pi})$. Hence, the final version of the problem of the sender is

\begin{equation}
\label{eq:finalformulation}
    \max_{\pi\in\mathrom{P}_{\mathcal{S}}(\mathcal{M})}\left(\sum_{(m,t,s)\in\mathcal{M}\times\mathcal{T}\times S}v(s,a_t^\star(\nu_{m,\pi}))\nu_{m,\pi}(s) \eta(t)p(m) \right).
\end{equation}

\subsection{Limitations of the model}
In this section, we discuss some limitations of the model.

The first caveat is somewhat technical: we have not explicitly emphasized that the mapping $(\nu, t) \mapsto a_t^\star(\nu)$ must be measurable. However, this requirement is essential for the analysis to hold. Although it is theoretically possible to construct non-measurable mappings of the form $(\nu, t) \mapsto a_t^\star(\nu)$, insisting on such pathological choices would be an unnecessarily adversarial stance.

A second consideration concerns the rule governing receivers when they are indifferent amongst several actions. The convention adopted here is that 'the receiver pleases the sender.' This choice has the advantage of making the sender's utility upper semicontinuous with respect to the communication policy (see Lemma \ref{lem:Wusc} below). An alternative rule is that, in cases of indifference, the receiver selects an action according to a predetermined distribution over the action set. This latter approach is closely related to the notion of regularization discussed in this article.

A third one is that the decision-making framework described in the previous section is not yet fully specified, as it may admit seemingly inconsistent action choices. This issue has significant implications, since the total number of admissible actions plays a critical role in estimating the optimal number of messages (see Theorem~\ref{thm:tausupp} below). In our setting, the space of all possible action profiles is given by
\[
\mathcal A^\mathcal T = \bigotimes_{t\in \mathcal T} \mathcal A = \left\{ (a_t)_{t\in \mathcal T} \text{ such that } a_t \in \mathcal A, \; \forall t \in \mathcal T \right\},
\]
which has cardinality $A^T$. However, allowing such a vast array of possible action profiles seems unnecessarily permissive. Indeed, suppose we have two different types $t_1$ and $t_2$, and a group of receivers $\mathcal{R} \subset \mathcal{L}$ that share the same  type, that is $(t_1)_{\mathcal{R}} = (t_2)_{\mathcal{R}}$. Further assume that the actions of all other receivers are identical across the two types, i.e., $\left(a^\star_{t_1}(\nu)\right)_{-\mathcal{R}} = \left(a^\star_{t_2}(\nu)\right)_{-\mathcal{R}}$. Then it seems reasonable?both from a modeling and intuitive standpoint?that the receivers in $\mathcal{R}$ should take the same action in both cases. Formally, we require that for all $t_1, t_2$ and all $\mathcal{R} \subset \mathcal{L}$,
\begin{equation}
\label{eq:restrictavailactions}
\left((t_1)_{\mathcal{R}}=(t_2)_{\mathcal{R}} 
\text{ and } \left(a^\star_{t_1}(\nu)\right)_{-\mathcal{R}}=\left(a^\star_{t_2}(\nu)\right)_{-\mathcal{R}}
\right)\Rightarrow 
\left(a^\star_{t_1}(\nu)\right)_{\mathcal{R}}=\left(a^\star_{t_2}(\nu)\right)_{\mathcal{R}},
\end{equation}
where, by $a_{-\mathcal{R}}$, we describe the actions taken by every receiver except for the ones in $\mathcal{R}$. Note, however, that it is entirely possible to fail to verify condition~\eqref{eq:restrictavailactions}, as illustrated in Section \ref{sec:exlimitactions}. As stated in said section, a way to circumvent this issue is to suppose that the sender is never indifferent to the actions of the receivers.

\subsection{The concavification of the sender's problem and the revelation principle}
In this section, following the ideas of \cite{kamenica2011bayesian}, we reformulate the sender's problem~\eqref{eq:finalformulation} into a linear programming problem set on the space of measures. For that purpose, for any prior $\nu \in \mathcal P(S)$, denote $W(\nu)$ the gain of the sender defined as
\begin{equation}
\label{eq:defiW}
    W(\nu) = \left(\sum_{(s,t)\in\mathcal{S}\times\mathcal{T}}v(s,a_t^\star(\nu)) \nu(s) \eta(t) \right).
\end{equation}
Then  the sender's problem~\eqref{eq:finalformulation} is reformulated in 
\begin{equation}
    \max_{\pi\in\mathrom{P}_{\mathcal{S}}(\mathcal{M})}\left(\sum_{m\in\mathcal{M}}W(\nu_{m,\pi}) p(m) \right).
\end{equation}
The first trick in the concavification of the problem relies on the introduction of $\tau_\pi$ which is a sum of Dirac masses at points $\nu_{m,\pi}$, that is
\[\tau_\pi=\sum_m p_m \delta_{\nu_{m,\pi}},\]
and to recast the problem into
\begin{equation}
    \max_{\pi\in\mathrom{P}_{\mathcal{S}}(\mathcal{M})}\int_{\mathrom {P}(\mathcal S)}  W(\nu) d \tau_\pi(\nu).
\end{equation}
The question is to describe the set of admissible $\tau_\pi$ when $\pi$ spans $\mathrom{P}_{\mathcal{S}}(\mathcal{M})$. First remark that 
because  $\sum p(m)=1$, then $\tau_\pi \in \mathrom {P}(\mathrom {P}(\mathcal S))$. Moreover , the equation $\sum_m \nu_{m,\pi}(s)=\mu(s)$ ensures that $\int_{\mathrom {P}(\mathcal S)} \nu d\tau_{\pi}(\nu)=\mu$.
We are led to introduce $\mathrom{T}_{M,\mu}$, the set of probability measure on $\mathrom {P}(\mathcal S)$ with expectation given by $\mu$ and convex combination of at most $M$ Dirac masses.
\begin{equation*}
    \mathrom{T}_{M,\mu}=\left\{\tau\in \mathrom P({\mathrom P (S)}) \text{ such that } \int_{\mathrom{P}(\mathcal{S})} \nu d\tau (\nu)=\mu\text{ and }  \#(\textrm{supp}(d\tau)) \le M \right\}.
\end{equation*}
For any $\pi \in \mathrom{P}_{\mathcal{S}}(\mathcal{M})$, then $\tau_{\pi}$ belongs to $\mathrom{T}_{M,\mu}$.  Reciprocally, if $\tilde \tau=\sum_{m}\rho(m)\delta_{\nu_m}$ is chosen in $\mathrom{T}_{M,\mu}$, then constructing $\pi(m,s)=\rho(m) \nu_m(s)/\mu(s)$, it is easy to check that $\pi\in \mathrom{P}_{\mathcal{S}}(\mathcal{M})$ and $\tilde \tau=\tau_\pi$. Hence, by a change of variable, the sender problem is equivalent to the  linear programming problem
\begin{equation}
\label{eq:senderprob:tausupp} 
    \max_{\tau \in\mathrom{T}_{M,\mu}}\int_{\mathrom {P}(\mathcal S)}  W(\nu) d \tau(\nu),
\end{equation}
A companion problem to \eqref{eq:senderprob:tausupp} is the one where the constraint on the number of messages is relaxed. We are then lead to consider the set $\mathrom{T}_{\mu}$, which is the set of probability measure on $\mathrom {P}(\mathcal S)$ with expectation given by $\mu$ 
\begin{equation}
    \mathrom{T}_{\mu}=\left\{\tau\in \mathrom P({\mathrom P (S)}) \text{ such that } \int_{\mathrom{P}(\mathcal{S})} \nu d\tau (\nu)=\mu \right\}.
\end{equation}
and the problem 
\begin{equation}
\label{eq:senderprob:tau}
    \hat W(\mu)=\max_{\tau \in\mathrom{T}_{\mu}}\int_{\mathrom {P}(\mathcal S)}  W(\nu) d \tau(\nu).
\end{equation}
Foreshadowing the developments in Section~\ref{sec:results} and Theorem~\ref{thm:tausupp}, Problem~\eqref{eq:senderprob:tau} is, in fact, equivalent to Problem~\eqref{eq:senderprob:tausupp} when the parameter $M$ is sufficiently large. In this regime, a version of the revelation principle applies and each message is a prescribed action. Problem~\eqref{eq:senderprob:tau} turns into the following linear programming problem :
\begin{equation}
\label{pbm:fulllinear}
\max_{
\pi \in \mathrom P_S(\mathcal A^\mathcal T), \pi \text{ verifies } \eqref{comm:pi}
} \sum_{(a,s) \in \mathcal A^\mathcal T\times \mathcal S} \sum_{t\in \mathcal T}v(s,a_t)\pi(a|s)\mu(s)\eta(t),
\end{equation}
where \eqref{comm:pi} is defined as
\begin{equation}
\label{comm:pi}
\sum_{s}u_{t_\ell}(s,(a_{t})_\ell)\pi(a|s)\mu(s) \ge 
\sum_{s}u_\ell(s,\tilde a)\pi(a|s)\mu(s),
\end{equation}
for any  $t\in \mathcal T$, $\ell \in\mathcal L$, $\tilde a \in  \mathcal A_\ell$ and $a\in \mathcal A^\mathcal T$.
\subsection{Regularization by strict convexity}
\subsubsection{Problems with the standard formulation}

In this section, we examine several shortcomings of the standard formulations \eqref{eq:senderprob:tau} and \eqref{pbm:fulllinear}, which motivate our proposal to regularize the receiver's problem. Although Problem~\eqref{eq:senderprob:tau} is linear in form, it is defined over a space of probability measures. Recall that for any function $f: A \to \mathbb{R}$, maximizing $f$ is equivalent to maximizing $\int f \, d\mu$ over $\mu \in \mathbb{P}(A)$. Consequently, linear programming problems over $\mathbb{P}(A)$ are, in general, as challenging as directly maximizing a non-linear function over $A$. In our setting, since the set $\mathcal{S}$ is finite, the space $\mathbb{P}(\mathcal{S})$ can be identified with the standard simplex in dimension $S+1$, and Problem~\eqref{eq:senderprob:tau} is therefore as difficult as maximizing a function over a space of dimension $S+1$.

For formulation \eqref{pbm:fulllinear}, the problem is linear in finite dimension. However, it involves $A^T \times S$ variables and approximately $A^{T+1}$ constraints. Consequently, the sheer scale of variables and constraints often renders it unsolvable in practice. In addition, if the solution set happens to be a facet, the algorithm becomes highly sensitive to parameter changes, as even minimal variations make it switch between vertices, which destabilizes the problem-solving process.

\subsubsection{Regularization}

In the previous model, the receiver acts as a rational agent in choosing its action with respect to its utility. We now introduce a variant where he is biased in his choice of strategy and wishes to remain close to an arbitrary strategy (the "irrational strategy"). To do so, we rely on a divergence function on the space $\mathsf{P}(\mathcal{S})$ which will quantifies the dissimilarity between two elements of $\mathsf{P}(\mathcal{S})$. The $\ell$-th receiver problem \eqref{eq:probinitreceiver:argmax} is reformulated as follows
\begin{equation}\label{eq:probregreceiver:util}
    \underset{\theta\in \mathsf{P}(\mathcal A_\ell)}{\mathrm{max}}\;\sum_{(s,a)\in\mathcal{S}\times\mathcal{A}_\ell} u_{t_{\ell}}(s,a)\theta(a)\nu(s) - \varepsilon \vartheta(\theta,\lambda_{\ell}),
\end{equation}
for some "irrational strategy" $\lambda_{\ell}\in\mathsf{P}(\mathcal A_\ell)$ and $\varepsilon>0$. In practice, we choose $\vartheta$ to be the Kullback-Leibler divergence given by
\begin{equation*}
    \vartheta(\nu,\mu) = \left\{\begin{array}{ll}
    \sum_{s\in\mathcal{S}}\nu(s)\log\left(\frac{\nu(s)}{\mu(s)}\right),\:\mbox{if }\nu\ll\mu,
    \\\infty,\:\mbox{else}.\end{array}\right.
\end{equation*}
The parameter $\varepsilon$ quantifies the degree of commitment to the irrational strategy.
\begin{remark}
A related model can be found in \cite{matyskova2023bayesian} where the authors introduce a costly information acquisition for the receiver. The cost is given by Shannon's entropy and connects to our model.
\end{remark}

We assume that $\mathrm{supp}(\lambda_{\ell}) = \mathcal{A}_{\ell}$, for any $\ell\in\{1,\ldots,L\}$. The previous problem is strictly concave and admits a unique solution denoted $\theta^{\star,\varepsilon}_{t_{\ell},\nu^{\pi}_m}$. We denote, for any $(a,t)\in\mathcal{A}\times\mathcal{T}^L$ and any $\nu\in\mathsf{P}(\mathcal{S})$,
\begin{equation}\label{eq:defthetareg0}
    \theta^{\star,\varepsilon}_{t,\nu}(a) = \prod_{\ell = 1}^L\theta^{\star,\varepsilon}_{t_{\ell},\nu}(a_{\ell}).
\end{equation}
By following the same arguments as in the previous section and by denoting, for any $a\in \mathcal{A}^{\mathcal{T}}$ and any $\nu\in\mathsf{P}(\mathcal{S})$,
\begin{equation}\label{eq:defthetareg}
    \theta_{\nu}^{\star,\varepsilon}(a) = \sum_{t\in\mathcal{T}}\theta_{t,\nu}^{\star,\varepsilon}(a)\eta(t)\quad\text{and}\quad W^{\varepsilon}(\nu) = \sum_{(s,a)\in\mathcal{S}\times\mathcal{A}^{\mathcal{T}}} w(s,a)\nu(s)\theta^{\star,\varepsilon}_{\nu}(a),
\end{equation}
the sender's problem becomes
\begin{equation}\label{eq:senderprobreg:tau}
    \max_{\tau\in \mathsf{T}_{\mu}}\left(\int_{\mathsf{P}(\mathcal{S})} W^{\varepsilon}(\nu)d\tau(\nu) \right).
\end{equation}

In contrast to Problem \eqref{eq:senderprob:tau}, this problem cannot be reformulated as a linear programming problem, since we have $\mathrm{supp}(\theta^{\star,\varepsilon}_{t,\nu}) = \mathcal{A}$. Nevertheless, it is possible to derive an explicit expression for $\theta^{\star,\varepsilon}_{t,\nu}$, which obviously depends on the choice of the divergence $\vartheta$. As we will see later in Section \ref{sec:Numericalcon}, first- or second-order methods can then be employed to solve this nonlinear optimization problem. In particular, this circumvent some of the limitations described in the previous section.

\section{Main results}\label{sec:results}
\subsection{On the concavification}

We now provide some results concerning Problem \eqref{eq:senderprob:tausupp} and Problem \eqref{eq:senderprob:tau}. To begin with, their constraint spaces are compact as stated in the following lemma.
\begin{lemma}\label{lem:compactau}
Let $M\in\mathbb{N}^*$. The spaces $\mathsf{T}_{\mu}$ and $\mathsf{T}_{M,\mu}$ are compact.
\end{lemma}

To address Problem \eqref{eq:senderprob:tau}, the role of concavification becomes clear through the following result, which both establishes the existence of a solution and provides a method for computing the optimal value through bi-conjugate Fenchel transform of $W$.

\begin{lemma}\label{lem:Wusc}
The function $W$ defined by \eqref{eq:defiW} is upper semicontinuous. The function $\hat W$ is the smallest concave function which is greater or equal to $W$. It is upper semicontinuous and is then equal to the bi-conjugate Fenchel transform of $W$.
\end{lemma}

A natural question when studying \eqref{eq:senderprob:tau} and \eqref{eq:senderprob:tausupp} is to determine the minimal value of $M$ (the number of messages) such that a solution to \eqref{eq:senderprob:tausupp} is also a solution to \eqref{eq:senderprob:tau}. 
\begin{remark}\label{rmk:representerthm}
Problem \eqref{eq:senderprob:tau} is a linear programming problem with $S$ constraints. As such, we know that it admits a solution $\tau\in\mathfrak{S}_{S}$ \cite{boyer2019representer}.
\end{remark}
The path to answering this question begins with the following decomposition lemma (Lemma~\ref{lem:decompotau}) concerning optimal solutions and its corollary (Corollary~\ref{cor:concavWtau}). 
\begin{lemma}\label{lem:decompotau}
Let $\tau \in \mathrom{T}_\mu$ be a solution to \eqref{eq:senderprob:tau} such that there exist $r \in [0,1]$, $\tau_1 \in \mathrom{T}_{\mu_1}$, and $\tau_2 \in \mathrom{T}_{\mu_2}$ satisfying $\tau = r\tau_1 + (1 - r)\tau_2$, so that $\mu = r\mu_1 + (1-r)\mu_2$. Then, $\tau_1$ and $\tau_2$ are also solutions to \eqref{eq:senderprob:tau} within their respective constraint sets.
\end{lemma}
\begin{corollary}\label{cor:concavWtau}
Let $\tau=\sum_{m=1}^M \rho_m \delta_{\nu_m}\in\mathrom{T}_{M,\mu}$ be a solution of \eqref{eq:senderprob:tau}. Then, for each $m$, we must have $W(\nu_m) = \hat{W}(\nu_m)$.
\end{corollary}
Given the latest corollary, we are in position to answer to the question of the optimal number of messages with a sharp estimate

\begin{theorem}\label{thm:tausupp}
There exists a solution $\tau^{\star}$ to \eqref{eq:senderprob:tau} in $\mathrom{T}_{\mu}$ that is a finite sum of Dirac measures supported on points $(\nu_i)_i \in \mathrom{P}(\mathcal{S})$, such that for any two distinct indices $i \ne j$, the induced actions $a^\star_t(\nu_i)$ and $a^\star_t(\nu_j)$ differ for at least one type of one receiver. Since there are at most $A^T$ such action profiles, problems \eqref{eq:senderprob:tau} and \eqref{eq:senderprob:tausupp} are equivalent whenever $M \ge A^T$. Conversely, for any given sets of receivers $\mathcal L$ and actions $\mathcal{A}$, if there is only one type, there exist sets of states $\mathcal{S}$ and utility functions $u$ and $v$ such that no solution to \eqref{eq:senderprob:tau} exists with support of cardinality less than or equal to $A$.
\end{theorem}

Theorem~\ref{thm:tausupp} is a revelation principle that states that there is a solution where each message $m=(a_t)_{t\in \mathcal T}$ is an element of $\mathcal A^\mathcal T$ of the form "If the configuration of the types is $t$, then perform the action $a_t \in \mathcal A$". Of course this message has to be credible, that is, once the posterior of this communication policy is computed, the action $a_t$ has to be admissible for every receiver, that is the communication policy $\pi$ has to verify, for every $a\in \mathcal A^\mathcal T$
\begin{equation}
\label{comm:pibis}
\sum_{s}u_{t_\ell}(s,(a_{t})_\ell)\pi(a|s)\mu(s) \ge 
\sum_{s}u_{t_\ell}(s,\tilde a)\pi(a|s)\mu(s), \forall t,\ell \in \mathcal T\times \mathcal L, \forall \tilde a \in  \mathcal A_\ell.
\end{equation}
We can restrict ourselves to studying such communication policies, and we obtain the following linear programming problem :
\begin{equation}
\label{pbm:fulllinearbis}
\max_{
\pi \in \mathrom P_S(\mathcal A^\mathcal T), \pi \text{ verifies } \eqref{comm:pibis}
} \sum_{(a,s) \in \mathcal A^\mathcal T\times \mathcal S} \sum_{t\in \mathcal T}v(s,a_t)\pi(a|s)\mu(s)\eta(t).
\end{equation}

\subsection{On the regularization}

We now turn to the regularized problem.

\begin{lemma}\label{lem:limsupvaluefun}
Let $(\tau^{\varepsilon})_{\varepsilon>0}$ be a sequence in $\mathsf{T}_{M,\mu}$, with $M\in\mathbb{N}$, that converges to some $\tau^0\in\mathsf{T}_{M,\mu}$.
Then, we have
\begin{equation*}
    \underset{\varepsilon\to0}{\mathrm{limsup}}\;\int_{\mathsf{P}(\mathcal{S})}W^{\varepsilon}(\nu) d\tau^{\varepsilon}(\nu)\leq \int_{\mathsf{P}(\mathcal{S})}W(\nu) d\tau^{0}(\nu).
\end{equation*}
\end{lemma}

By using the previous lemma with a sequence of maximizers (which have a finite support of cardinality $S$, by Remark \ref{rmk:representerthm}), we directly deduce the following corollary.

\begin{corollary}\label{cor:suboptregtau}
For any $\varepsilon>0$, denote 
\begin{equation*}
    \bar{W}^{\varepsilon,\star} = \max_{\tau\in \mathsf{T}_{\mu}}\left(\int_{\mathsf{P}(\mathcal{S})} W^{\varepsilon}(\nu)d\tau(\nu) \right) \quad\text{and}\quad \bar{W}^{\star} = \max_{\tau\in \mathsf{T}_{\mu}}\left(\int_{\mathsf{P}(\mathcal{S})} W(\nu)d\tau(\nu) \right).
\end{equation*}
Then, we have
\begin{equation*}
    \underset{\varepsilon\to0}{\mathrm{limsup}}\;\bar{W}^{\varepsilon,\star}\leq \bar{W}^{\star}.
\end{equation*}
\end{corollary}

In order to obtain the convergence of $\bar{W}^{\varepsilon,\star}$ to $\bar{W}^{\star}$ as $\varepsilon$ goes to $0$, we require the following assumption.

\begin{assumption}\label{asm:reg}
If $a\in \mathcal A^\mathcal T$ is a prescribed action for a certain prior, that is there exists $\nu\in\mathsf{P}(\mathcal{S})$ such that $a_t=a^\star_t(\nu)$, for any $t\in\mathcal{T}$, then it is a forced action for a perhaps different prior with a support contained in $\mathrm{supp}(\mu)$. That is, there exists $\tilde \nu\in\mathsf{P}(\mathcal{S})$ such that $\mathrm{supp}(\tilde \nu)\subset\mathrm{supp}(\mu)$ and
$\mathcal A_t(\tilde \nu)=\{a_t\}$, for any $t\in\mathcal{T}$.
\end{assumption}

\begin{lemma}\label{lem:convergtaureg}
Under Assumption \ref{asm:reg}, there exists a sequence $(\tau^{\varepsilon})_{\varepsilon>0}\subset \mathsf{T}_{\mu}$ such that
\begin{equation*}
    \int_{\mathsf{P}(\mathcal{S})}W^{\varepsilon}(\nu) d\tau^{\varepsilon}(\nu)\underset{\varepsilon\to0}{\longrightarrow}\bar{W}^{\star}.
\end{equation*}
\end{lemma}

\begin{corollary}\label{cor:convergWeg}
Under Assumption \ref{asm:reg}, we have
\begin{equation*}
    \lim_{\varepsilon\to0}\bar{W}^{\varepsilon,\star} =\bar{W}^{\star}.
\end{equation*}
\end{corollary}

\begin{remark}
Assumption \ref{asm:reg} is essential for the previous result as shown in Section \ref{ex:SharpAsm} below.
\end{remark}

As a result, we have the following theorem.

\begin{theorem}\label{theo:conv:regul}
Under Assumption \ref{asm:reg}, let $(\tau^{\varepsilon})_{\varepsilon>0}$ be a sequence of maximizers of \eqref{eq:senderprobreg:tau} such that, for some $M\in\mathbb{N}$ and any $\varepsilon>0$, $\tau^{\varepsilon}\in\mathsf{T}_{M,\mu}$. Then, this sequence converges, up to a subsequence, to a maximizer of \eqref{eq:senderprob:tau} as $\varepsilon$ goes to $0$.
\end{theorem}

\subsection{Numerical considerations}\label{sec:Numericalcon}
In this section, we provide a rationale for some of the numerical choices underlying our method. We begin with the now-classic $\softmax$ function, which takes as input any vector $x \in \mathbb{R}^n$ and returns an element of the probability simplex $\Delta_n$, defined componentwise as
\[
\softmax(x)_i = \frac{e^{x_i}}{\displaystyle\sum_{j=1}^n e^{x_j}} \quad \text{for all } 1 \le i \le n.
\]
This function is a staple in neural network classification tasks, typically serving as the final output layer. The reasons are straightforward: it provides a smooth, numerically stable, and easily differentiable way to turn scores into probabilities. Moreover, it plays exceptionally well with the Kullback-Leibler divergence, making it a natural fit in probabilistic modeling.

To be more precise, if $\vartheta$ is chosen as  the Kullback-Leibler divergence then the explicit solution of the optimization problem \eqref{eq:probregreceiver:util} is given by

\begin{equation}
\label{eq:thetabysoftmax}
\theta^{\star,\varepsilon}_{t_\ell,\nu}(\cdot)=\lambda_\ell(\cdot)\, \softmax\left(\frac{\sum_s u_{t_\ell}(s,\cdot)\nu(s)}{\varepsilon}\right)
\end{equation}
\paragraph{Choice of optimization variable.}
In our implementation, we invoke the $\softmax$ function not once, but twice. Instead of directly optimizing over $\pi$ under the usual probabilistic constraints, we adopt a relaxed formulation: we assume $\pi = \softmax(x)$ for some unconstrained variable $x$, and we perform second-order optimization directly on $x$. This sidesteps the need for constrained optimization, which can often be more cumbersome.
In summary, the optimization pipeline unfolds as follows:
\begin{itemize}
    \item Given a function $x:\mathcal M\times \mathcal S\rightarrow \mathbb{R}$, compute $\pi \in \mathrom{P}_{\mathcal S}(\mathcal M)$ given by
    $$\pi(\cdot \mid s) =\softmax(x(\cdot,s)), \quad \forall s\in \mathcal S.$$
    \item For each message $m\in \mathcal M$, compute $\nu_{m,\pi}$ and $p(m)$ by \eqref{eq:define:nu} and then for each receiver $\ell\in \mathcal L$ of type $t_\ell\in \mathcal T_\ell$, compute $\theta^{\star,\varepsilon}_{t_\ell,\nu_{m,\pi}}$ by \eqref{eq:thetabysoftmax} then $\theta^{\star,\varepsilon}_{\nu_{m,\pi}}$ by \eqref{eq:defthetareg}.
    \item Finally, evaluate the objective 
    $${\bf v}(x)=\sum_{m\in\mathcal{M}}\sum_{(s,a)\in\mathcal{S}\times\mathcal{A}^{\mathcal{T}}} v(s,a)\nu_{m,\pi}(s)\theta^{\star,\varepsilon}_{\nu_{m,\pi}}(a)p(m).$$
    
\end{itemize}
To maximize the objective $x \mapsto {\bf v}(x)$, we rely on a standard BFGS algorithm with Wolfe line search?tried and tested tools for unconstrained smooth optimization. For gradient computations, we turn to PyTorch \cite{paszke2019pytorch}, which offers automatic differentiation and a highly optimized native implementation of the $\softmax$ function. This makes it an ideal companion for the task at hand: fast, reliable, and doing most of the hard work under the hood.

\paragraph{Choice of $\lambda_\ell$ and the regularization parameter $\varepsilon$.}
The weighting function $\lambda_\ell$ is quite flexible: any choice is admissible as long as $\lambda_\ell(a) > 0$ for all $a \in \mathcal{A}_\ell$. For simplicity?and to avoid introducing additional tuning parameters?we adopt the uniform distribution, i.e.,
\[
\lambda_\ell(a) = \frac{1}{A_\ell} \quad \text{for all } a \in \mathcal{A}_\ell.
\]

We now turn to the choice of the regularization parameter $\varepsilon$. A natural impulse is to take $\varepsilon$ as small as possible in order to better approximate the unregularized problem. However, there is a catch: when $\varepsilon$ becomes too small, the $\softmax$ function begins to closely approximate the $\arg\max$ operator. This is problematic, as the gradient of $\arg\max$ is either zero or undefined?both of which are highly undesirable in a gradient-based optimization routine. This effect is illustrated below in Section \ref{sec:hunterex}. In practice, if utility values are roughly of order one, a value of $\varepsilon = 10^{-2}$ provides a good compromise: it preserves enough smoothness for reliable optimization while remaining close enough to the original (non-regularized) objective to yield meaningful results. If small values of $\varepsilon$ are crucial, we propose in Section \ref{sec:hunterex} to iteratively reduce $\varepsilon$ using a strategy inspired by the interior-point method.

\section{Examples}\label{sec:Examples}
\subsection{Analytical examples}
\subsubsection{Sharpness of the assumption of Theorem \ref{theo:conv:regul}}\label{ex:SharpAsm}

This example aims to demonstrate a sequence of maximizers of the regularized problem whose utilities fail to converge to that of the original problem. In this setting, Assumption \ref{asm:reg} is not satisfied, and the conclusion of Corollary \ref{cor:convergWeg} does not hold. We consider $L = 1$, $\mathcal{S} = \{s_1,s_2\}$, $\mathcal{A}_1 = \{a_1,a_2\}$ and $\mathcal{T}_1 = \{t_1\}$. We set the divergence to be the Kullback-Leibler divergence with $\lambda = (1/2,1/2)$ and $\mu = (1/2,1/2)$. The utilities are
\begin{equation*}
    u = \begin{pmatrix}1 & 0 \\ 0 & 0 \end{pmatrix}\quad\mbox{and}\quad v = \begin{pmatrix}0 & 1 \\ 0 & 1 \end{pmatrix}.
\end{equation*}
In this case, we observe that, for any $\nu\in\mathsf{P}(\mathcal{S})$,
\begin{equation*}
    \sum_{s\in\mathcal{S}}u(s,a_1)\nu(s) = \nu(s_1)\quad\text{and}\quad \sum_{s\in\mathcal{S}}u(s,a_2)\nu(s) = 0.
\end{equation*}
Thus, we have $\mathcal{A}_{1,t_1}(\nu) = \{a_1\}$, for all $\nu\in\mathsf{P}(\mathcal{S})$ such that $\nu \neq \delta_{s_2}$, and $\mathcal{A}_{1,t_1}(\delta_{s_2}) = \{a_1,a_2\}$. Thus, Assumption \ref{asm:reg} does not hold. Indeed, there exists no $\nu\in\mathsf{P}(\mathcal{S})$ such that $\mathcal{A}_{1,t_1}(\nu) = \{a_2\}$. We have, for any $\nu\in\mathsf{P}(\mathcal{S})$,
\begin{equation*}
    \theta^{\star,\varepsilon}_{\nu} = \frac{e^{\nu(s_1)/\varepsilon}}{1 + e^{\nu(s_1)/\varepsilon}}\delta_{a_1} + \frac{1}{1 + e^{\nu(s_1)/\varepsilon}}\delta_{a_2},
\end{equation*}
and it follows that
\begin{equation*}
    \lim_{\varepsilon\to0}\theta^{\star,\varepsilon}_{\nu} = \left\{\begin{array}{cc}
         \delta_{a_1},&\text{ if }\nu\neq \delta_{s_2},\\
         (\delta_{a_1}+\delta_{a_2})/2,&\text{ if }\nu=\delta_{s_2}.
    \end{array}\right.
\end{equation*}
We can see that
\begin{equation*}
    W^{\varepsilon}(\nu) = \frac{1}{1 + e^{\nu(s_1)/\varepsilon}},
\end{equation*}
so that a solution of \eqref{eq:senderprobreg:tau} is
\begin{equation*}
    \tau = \frac12\delta_{\delta_{s_1}} + \frac12\delta_{\delta_{s_2}},
\end{equation*}
and it follows that $\bar{W}^{\star,\varepsilon} = \frac12\left(\frac{1}{1 + e^{1/\varepsilon}} + \frac{1}{2} \right)\to \frac14$ as $\varepsilon\to0$. However, we can see that
\begin{equation*}
    \theta^{\star}_{\nu} = \mathbbm{1}_{\{\nu(s_1)>0\}}\delta_{a_1} + \mathbbm{1}_{\{\nu(s_1)=0\}}\delta_{a_2},
\end{equation*}
which yields
\begin{equation*}
    W(\nu) = \mathbbm{1}_{\{\nu(s_1)=0\}}\quad\mbox{and}\quad \int_{\mathsf{P}(\mathcal{S})}W(\nu)d\tau(\nu) = 1/2 = \bar{W}^{\star}.
\end{equation*}
In the end, we have $\lim_{\varepsilon\to0}\bar{W}^{\star,\varepsilon} = 1/4<\bar{W}^{\star} = 1/2$.

\subsubsection{Inconsistency in admissible action profiles}\label{sec:exlimitactions}

%Consider two receivers, $a$ and $b$, and two types for receiver $b$, namely $b_1$ and $b_2$. Let there be two possible actions: $\alpha = (\alpha_1, \alpha_2)$ for receiver $a$, and $\beta = (\beta_1, \beta_2)$ for receiver $b$. Suppose that, under some peculiar prior $\nu$, both the sender and receiver $a$ are indifferent between the two components of $\alpha$. Now suppose that $b_1$ strictly prefers $\beta_1$, while $b_2$ remains indifferent between $\beta_1$ and $\beta_2$. Finally, assume that the sender also strictly prefers $\beta_1$.

%In the first type $(a, b_1)$, the set of admissible receiver actions is $\alpha \times \{\beta_1\}$. Since the sender is indifferent over $\alpha$, any pair in this set is optimal from her perspective. A decision must thus be made?say, arbitrarily, $(\alpha_1, \beta_1)$ is selected.

%In the second type $(a, b_2)$, the set of admissible actions is $\alpha \times \beta$, but among these, the sender prefers those in $\alpha \times \{\beta_1\}$. Once again, a choice must be made, and?perhaps just to mix things up?$(\alpha_2, \beta_1)$ is selected.

%This clearly violates condition~\eqref{eq:restrictavailactions}.

In this section, we demonstrate that the model may be inconsistent, in the sense that condition \eqref{eq:restrictavailactions} need not hold. We consider the case of two states, two receivers with the second one having two different types. That is, $L = 2$, $\mathcal{S} = \{s_1,s_2\}$, $\mathcal{A}_1 = \{a_{11},a_{12}\}$, $\mathcal{A}_2 = \{a_{21},a_{22}\}$, $\mathcal{T}_1 = \{t_{1}\}$ and $\mathcal{T}_2 = \{t_{21},t_{22}\}$. The utilities are the following
\begin{align*}
    &v(\cdot,a_{11},\cdot) = v(\cdot,a_{12},\cdot) = \begin{pmatrix} 1 & 0 \\ 1 & 0\end{pmatrix},\quad u_{t_1} = \begin{pmatrix} 0 & 0 \\ 0 & 0\end{pmatrix},\\ &u_{t_{21}} = \begin{pmatrix} 1 & 0 \\ 1 & 0\end{pmatrix}\quad\text{and}\quad u_{t_{22}} = \begin{pmatrix} 0 & 0 \\ 0 & 0\end{pmatrix}.
\end{align*}
For any state $s\in\mathcal{S}$, we observe that both the sender and the first receiver are indifferent to the choice of actions in $\mathcal{A}_1$. Furthermore, the sender strictly prefer action $a_{21}$. The second receiver of the first type also strictly prefer action action $a_{21}$ while the second type is completely indifferent.

Thus, for any $\nu\in\mathsf{P}(\mathcal{S})$, the admissible actions are
\begin{equation*}
    \mathcal{A}^{\star}_{t_1}(\nu) = \mathcal{A}_1,\quad \mathcal{A}^{\star}_{t_{21}}(\nu) = \{a_{21}\}\quad\text{and}\quad \mathcal{A}^{\star}_{t_{22}}(\nu) = \mathcal{A}_2.
\end{equation*}
If follows that there are several possible optimal actions in each case, but assume that the sender chooses
\begin{equation*}
    a^{\star}_{t_1,t_{21}}(\nu) = (a_{11},a_{21})\quad\text{and}\quad a^{\star}_{t_1,t_{22}}(\nu) = (a_{12},a_{21}),
\end{equation*}
which clearly violates condition~\eqref{eq:restrictavailactions}.

Now, consider again how the receivers arrive at their choices. As previously assumed, each receiver discloses their set of admissible actions, and together they select one that is most favorable to the sender. In this case, receiver $a$ realizes that he is facing two distinct situations (because the admissible actions of receiver $b$ differ), and so perceives no inconsistency. 

However, if the receivers privately communicate their admissible sets of actions to the sender, after which the sender selects and reveals a final action, receiver $a$ experiences an inconsistency because he observes $b$'s chosen action but not the set of actions available to him.

\subsubsection{Regularization of the judge and prosecutor example}

\begin{figure*}[h!]
    \centering
    \begin{subfigure}[t]{0.5\textwidth}
        \centering
        \includegraphics[width = 18em]{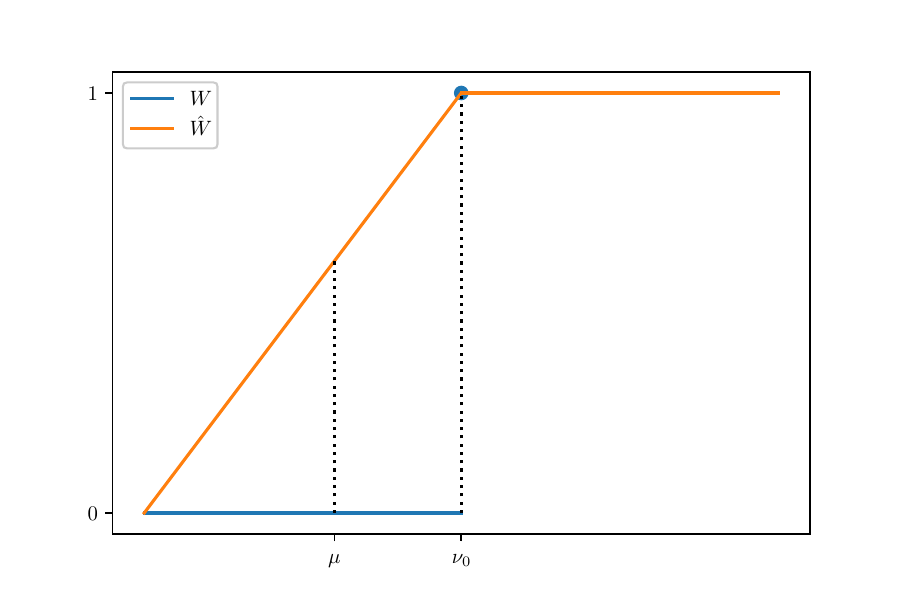}
        \caption{The unregularized problem}\label{fig:W}
    \end{subfigure}%
    ~ 
    \begin{subfigure}[t]{0.5\textwidth}
        \centering
        \includegraphics[width = 18em]{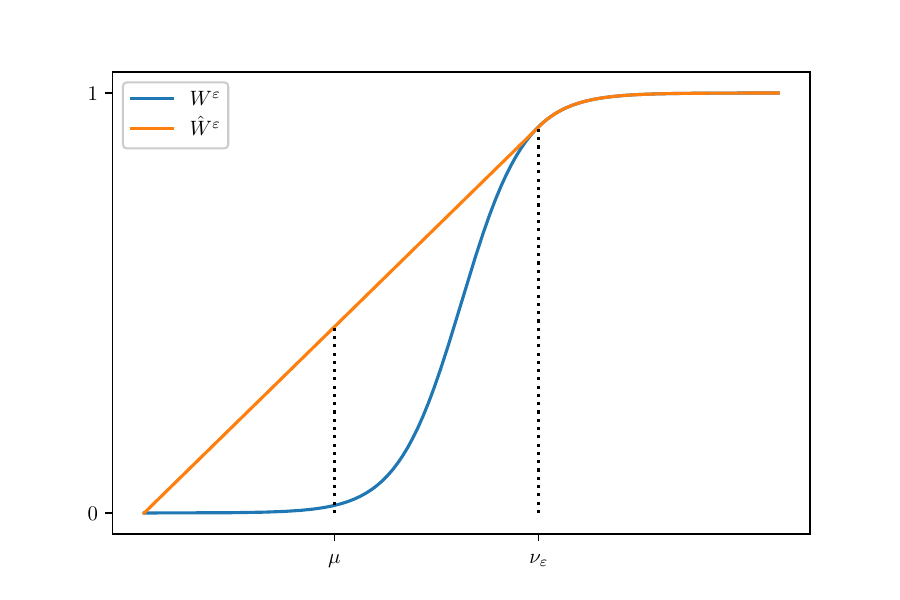}
        \caption{Regularization with a uniform irrational strategy $\lambda = (1/2,1/2)$}\label{fig:Wregunif}
    \end{subfigure}
    \\
    \centering
    \begin{subfigure}[t]{0.5\textwidth}
        \centering
        \includegraphics[width = 18em]{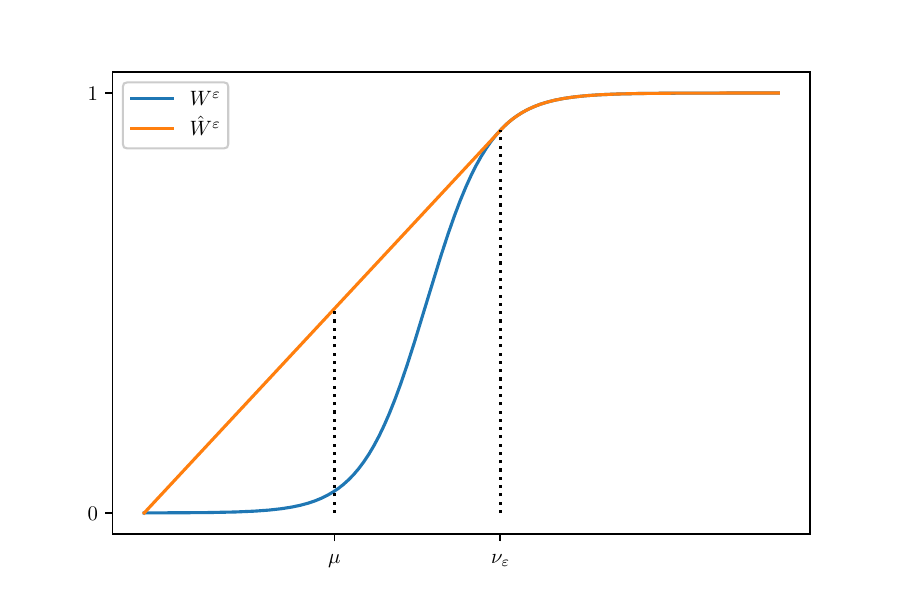}
        \caption{Regularization with a gullible irrational strategy $\lambda = (1/4,3/4)$}\label{fig:Wreggul}
    \end{subfigure}%
    ~ 
    \begin{subfigure}[t]{0.5\textwidth}
        \centering
        \includegraphics[width = 18em]{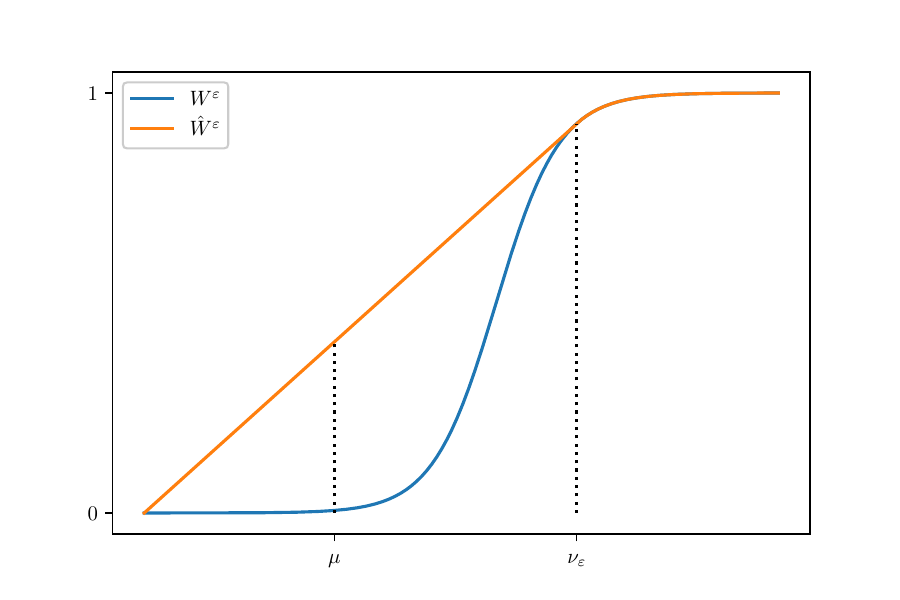}
        \caption{Regularization with a stubborn irrational strategy $\lambda = (3/4,1/4)$}\label{fig:Wregstub}
    \end{subfigure}
    \caption{The effects of regularization on the standard example of the judge and prosecutor for $\varepsilon = 0.1$ and different irrational strategies.}\label{fig:examplejudge}
\end{figure*}

In this example, we carry out a step-by-step analysis of the well-known prosecutor?judge example from \cite{kamenica2011bayesian}, applying a regularization based on the Kullback-Leibler divergence. We have a single receiver, one type, $\mathcal{S} =\{s_1,s_2\}$ and $\mathcal{A} =\{a_1,a_2\}$. The prior is $\mu = (0.7,0.3)$ and the utilities are given by
\begin{equation*}
    u = \begin{pmatrix}1 & 0 \\ 0 & 1 \end{pmatrix}\quad\mbox{and}\quad v = \begin{pmatrix}0 & 1 \\ 0 & 1 \end{pmatrix}.
\end{equation*}
We recall that, in this setting, the optimal solution to \eqref{eq:senderprob:tau} is given by
\begin{equation*}
    \tau^0 = (1-p_0)\delta_{\delta_{s_1}} + p_0\delta_{\nu_0},
\end{equation*}
where $p_0 = 2\mu(s_2)$ and $\nu_0 = (1/2,1/2)$. Furthermore, the functions $W$ and $\hat{W}$ are illustrated in Figure \ref{fig:W}.

For any $\nu\in\mathsf{P}(\mathcal{S})$, we compute
\begin{equation*}
    \sum_{s\in\mathcal{S}}u(s,a_1)\nu(s) = \nu(s_1) = 1 - \nu(s_2)\quad\mbox{and}\quad \sum_{s\in\mathcal{S}}u(s,a_2)\nu(s) = \nu(s_2),
\end{equation*}
so that
\begin{align*}
    \theta^{\star,\varepsilon}_{\nu} &= \frac{\lambda(a_1)}{\lambda(a_1) + \lambda(a_2)e^{(2\nu(s_2)-1)/\varepsilon}}\delta_{a_1} + \frac{\lambda(a_2)}{\lambda(a_2) + \lambda(a_1) e^{-(2\nu(s_2)-1)/\varepsilon}}\delta_{a_2} 
    \\ &= (1-\sigma_{\lambda,\varepsilon}(\nu(s_2)))\delta_{a_1} + \sigma_{\lambda,\varepsilon}(\nu(s_2))\delta_{a_2},
\end{align*}
where 
\begin{equation*}
    \sigma_{\lambda,\varepsilon}(x) = \frac{\lambda(a_2)}{\lambda(a_2) + \lambda(a_1) e^{-(2x-1)/\varepsilon}}.
\end{equation*}
Thus, we obtain
\begin{equation*}
    W^{\varepsilon}(\nu) = \sigma_{\lambda,\varepsilon}(\nu(s_2)).
\end{equation*}
The concavification of $W^{\varepsilon}$ is such that
\begin{equation*}
    \hat{W}^{\varepsilon}(\nu) =\left\{\begin{array}{cc}
          \sigma_{\lambda,\varepsilon}(r_{\varepsilon}) + \sigma_{\lambda,\varepsilon}'(r_{\varepsilon})(r_{\varepsilon} - \nu(s_2)),&\text{ if }\nu(s_2)\leq r_{\varepsilon}, \\
          \sigma_{\lambda,\varepsilon}(\nu(s_2)),&\text{ if }\nu(s_2)>r_{\varepsilon},
    \end{array}\right.
\end{equation*}
where $r_{\varepsilon}\in[0.5,1]$ is the solution of
\begin{equation*}
    \sigma_{\lambda,\varepsilon}(0) = \sigma_{\lambda,\varepsilon}(r_{\varepsilon}) + \sigma_{\lambda,\varepsilon}'(r_{\varepsilon})r_{\varepsilon}.
\end{equation*}
In particular, a solution of  \eqref{eq:senderprobreg:tau} is
\begin{equation*}
    \tau^{\varepsilon} =(1-p_{\varepsilon})\delta_{\delta_{s_1}} + p_{\varepsilon}\delta_{\nu_{\varepsilon}},
\end{equation*}
where $p_{\varepsilon} = \mu(s_2)/r_{\varepsilon}$ and $\nu_{\varepsilon} = (1-r_{\varepsilon},r_{\varepsilon})$. The function $W^{\varepsilon}$ and $\hat{W}^{\varepsilon}$, as well as $\nu_{\varepsilon}$, are depicted in Figures \ref{fig:Wregunif}-\ref{fig:Wreggul}-\ref{fig:Wregstub} for different irrational strategies $\lambda$. We observe that $\hat{W}^{\varepsilon}(\mu)\leq \hat{W}(\mu)$, for any choice of $\lambda$. Furthermore, $\hat{W}^{\varepsilon}(\mu)$ increases with $\lambda(a_2)$, which is naturally expected in this model.

\subsection{Numerical examples}

The numerical examples below are all done with \textit{BASIL} (Bayesian SIgnaling Library), a library in Python publicly avaiblable\footnote{https://plmlab.math.cnrs.fr/degourna/bayesian-persuasion-by-regularization}.

\subsubsection{The voting problem}\label{sec:votingproblem}
We first study the voting problem of \cite{alonso2016persuading} that consists in three receivers (voters A, voters B and voters C) which have to vote between two options. In this case the set of actions is the same for every receiver $\mathcal A=\{\text{Yes},\text{No}\}$, there are three states $\mathcal S = \{\text{A},\text{B},\text{C}\}$ and the utility $u_i(s,a)$ of the three receivers are given by

$$u_1=\begin{pmatrix} 1.1 & 0 \\ -1 &0 \\ -1&0 \end{pmatrix},
u_2=\begin{pmatrix} -1 & 0 \\ 1.1 &0 \\ -1&0 \end{pmatrix},
u_3=\begin{pmatrix} -1 & 0 \\ -1 &0 \\ 1.1&0 \end{pmatrix},
$$

\begin{figure}[htb!]
    \centering
    \includegraphics[width=0.8\textwidth]{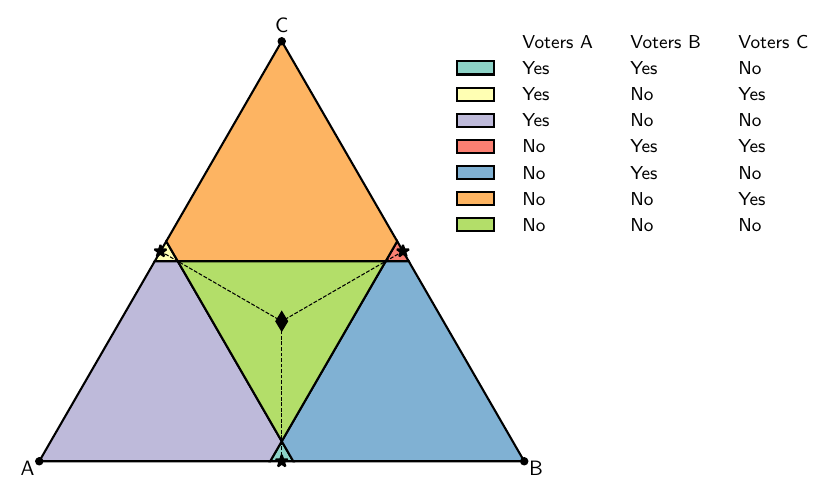}
    \caption{The voting problem. A point in the triangle correspond to an element of $\mathcal P(\mathcal S)$, the pure states are the vertices. The coloring corresponds to the different strategies of the receivers (voters) for a given state. The black diamond represents $\mu$ and the three stars represents the different $\nu$ for the three messages. Each star is located in one of the small triangle where the majority of voters approve the bill. In these triangles, the utility of the sender is $1$, it is $0$ elsewhere.}\label{fig:voters_actions}
\end{figure}

The prior is given by $\mu=(1/3,1/3,1/3)$. The actions taken by the three receivers, as well as the optimal signal, is displayed in Figure~\ref{fig:voters_actions}. 

In this Figure, the prior is a point in the simplex $\Delta_3$, where the pure states are the vertices of the triangle, beginning at the lower left and proceeding in counter-clockwise order.
The utility of the sender is $1$ if there is a majority of "Yes" and $0$ otherwise.  Without persuasion, each receiver will vote "No" and the corresponding gain of the sender is $0$. Consider the decomposition 
\[\mu=\frac 1 3 \sum_i \nu_i \quad \text{with } \nu_1=\left(0,\frac 12,\frac 12\right), \nu_2=\left(\frac 12,0,\frac 12\right) \text{ and } \nu_3=\left(\frac 12,\frac 12,0\right).\]
For each $\nu_i$ two receivers will vote "Yes", hence the sender gains $1$. For this particular choice of $d\tau= \sum_i \frac 1 3 \delta_{\nu_i}$, the sender reaches the maximum of the utility. This decomposition is represented with black stars in Figure~\ref{fig:voters_actions}.

\begin{figure}
    \centering
    \includegraphics[width=0.4\textwidth]{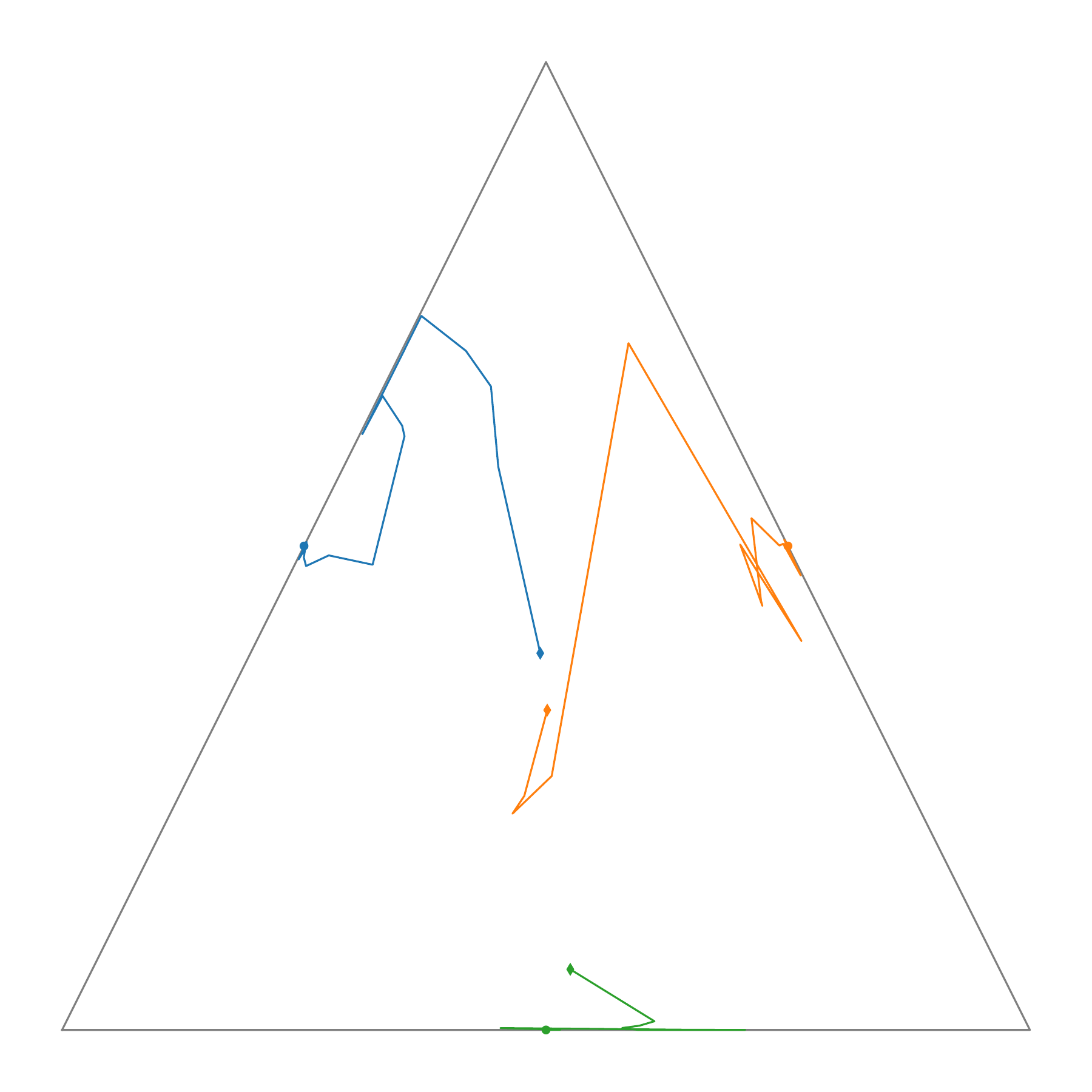}
    \includegraphics[width=0.4\textwidth]{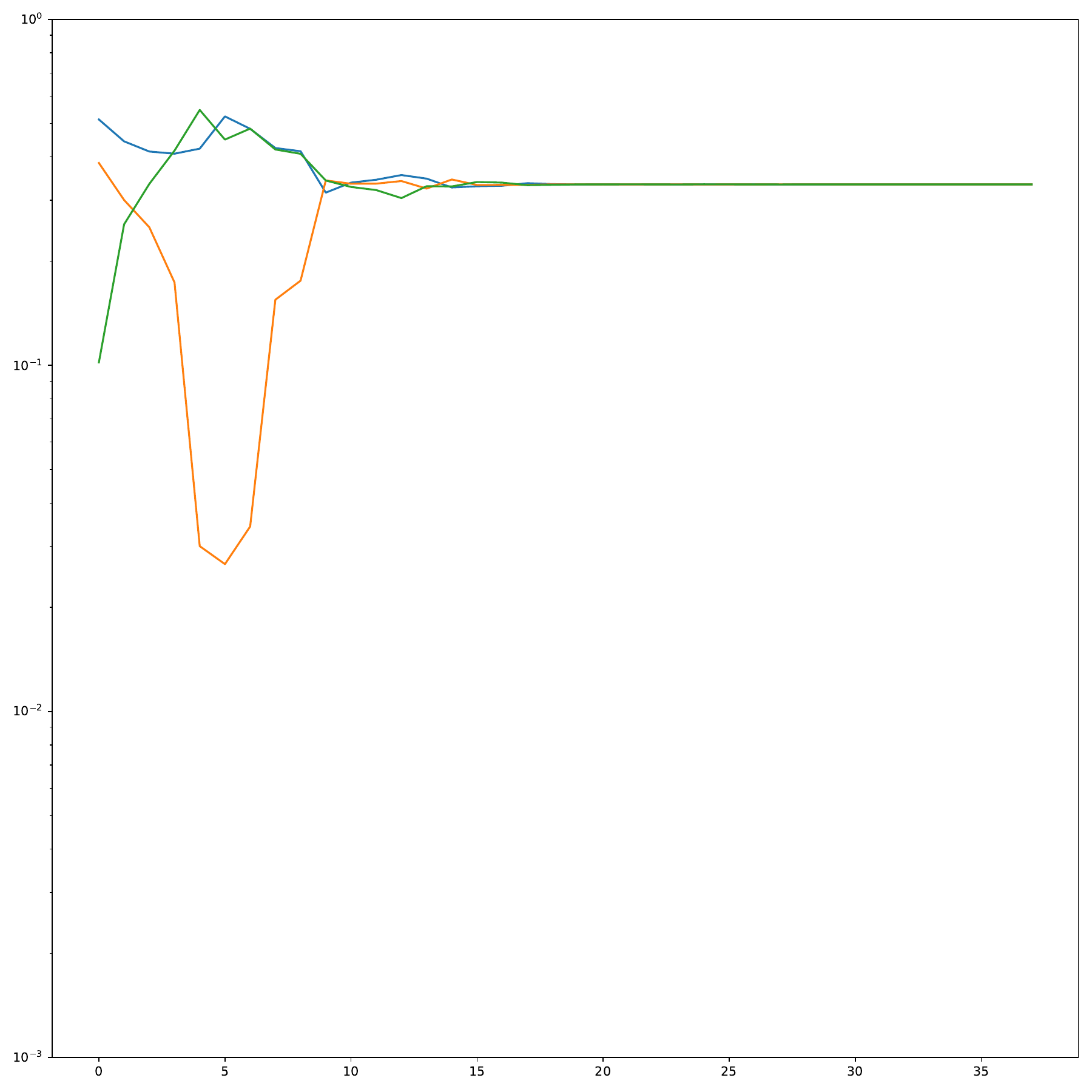}
    \includegraphics[width=0.4\textwidth]{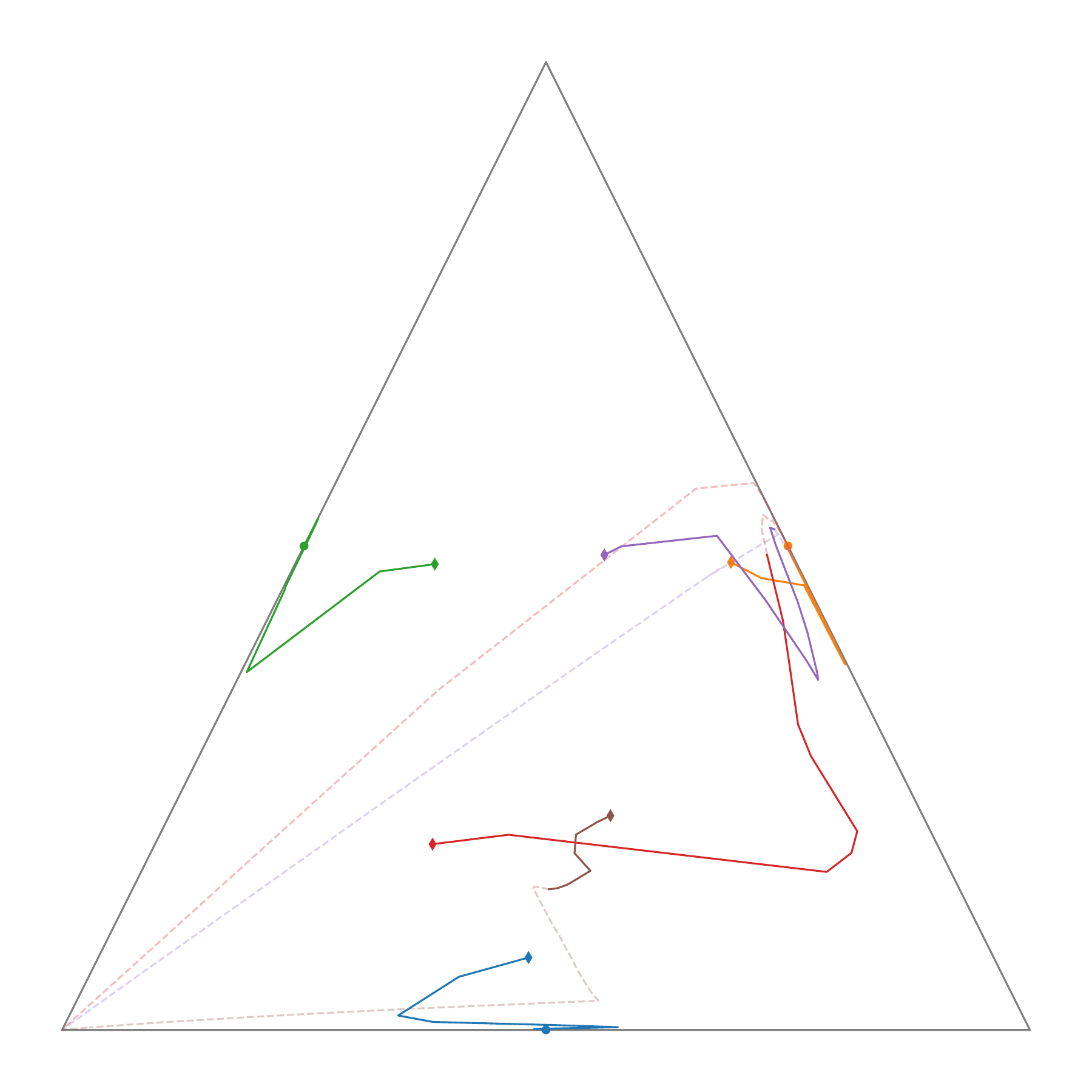}
    \includegraphics[width=0.4\textwidth]{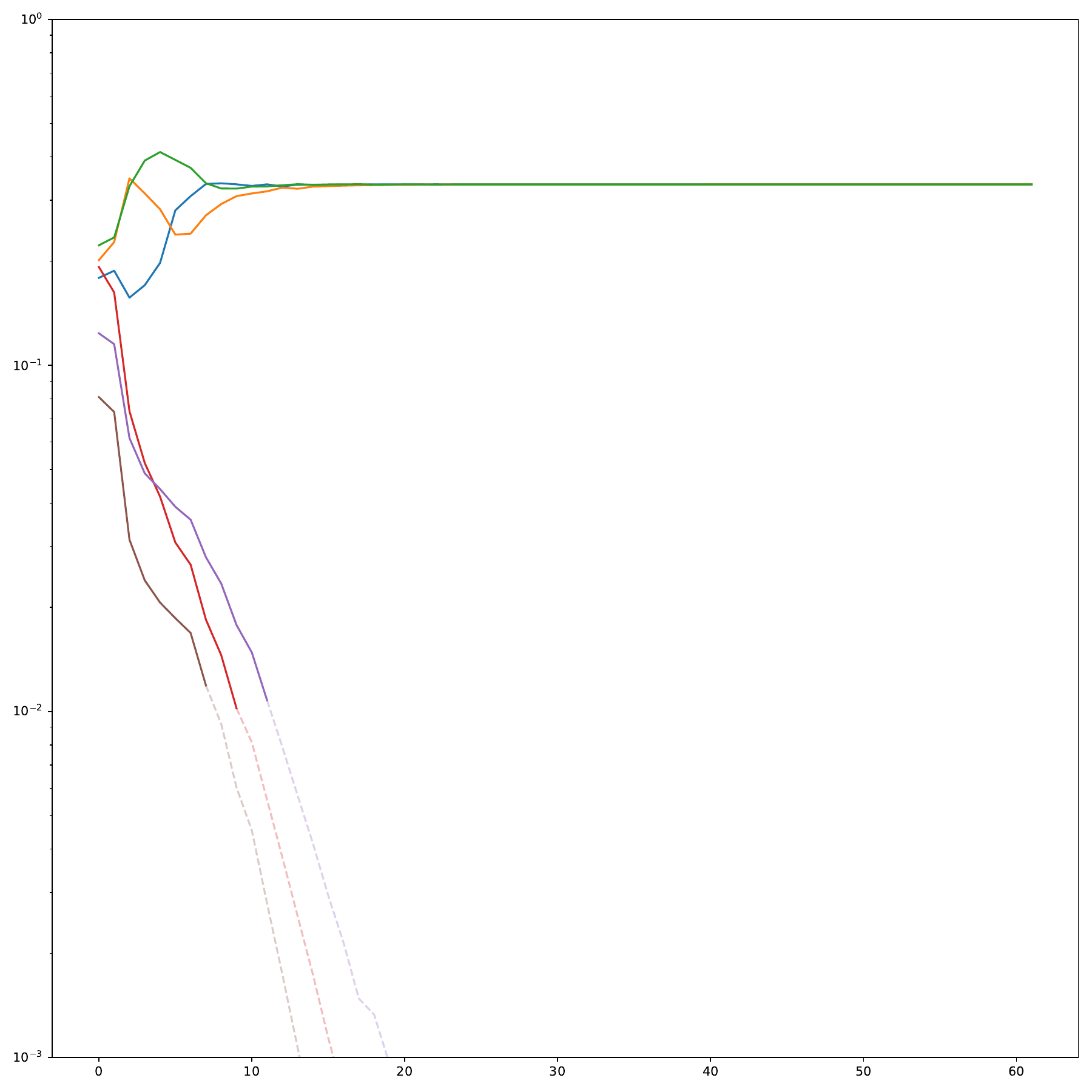}
    \includegraphics[width=0.4\textwidth]{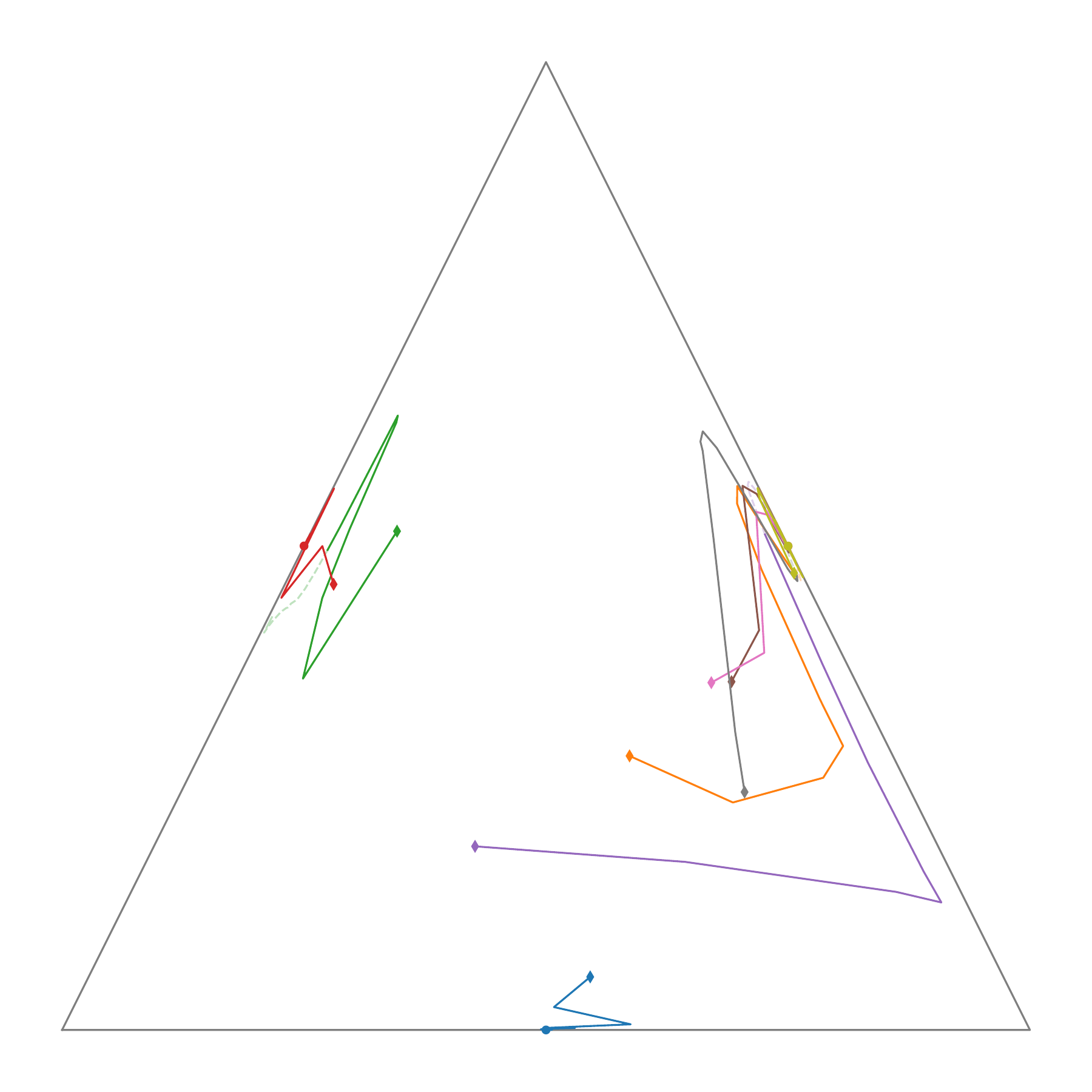}
    \includegraphics[width=0.4\textwidth]{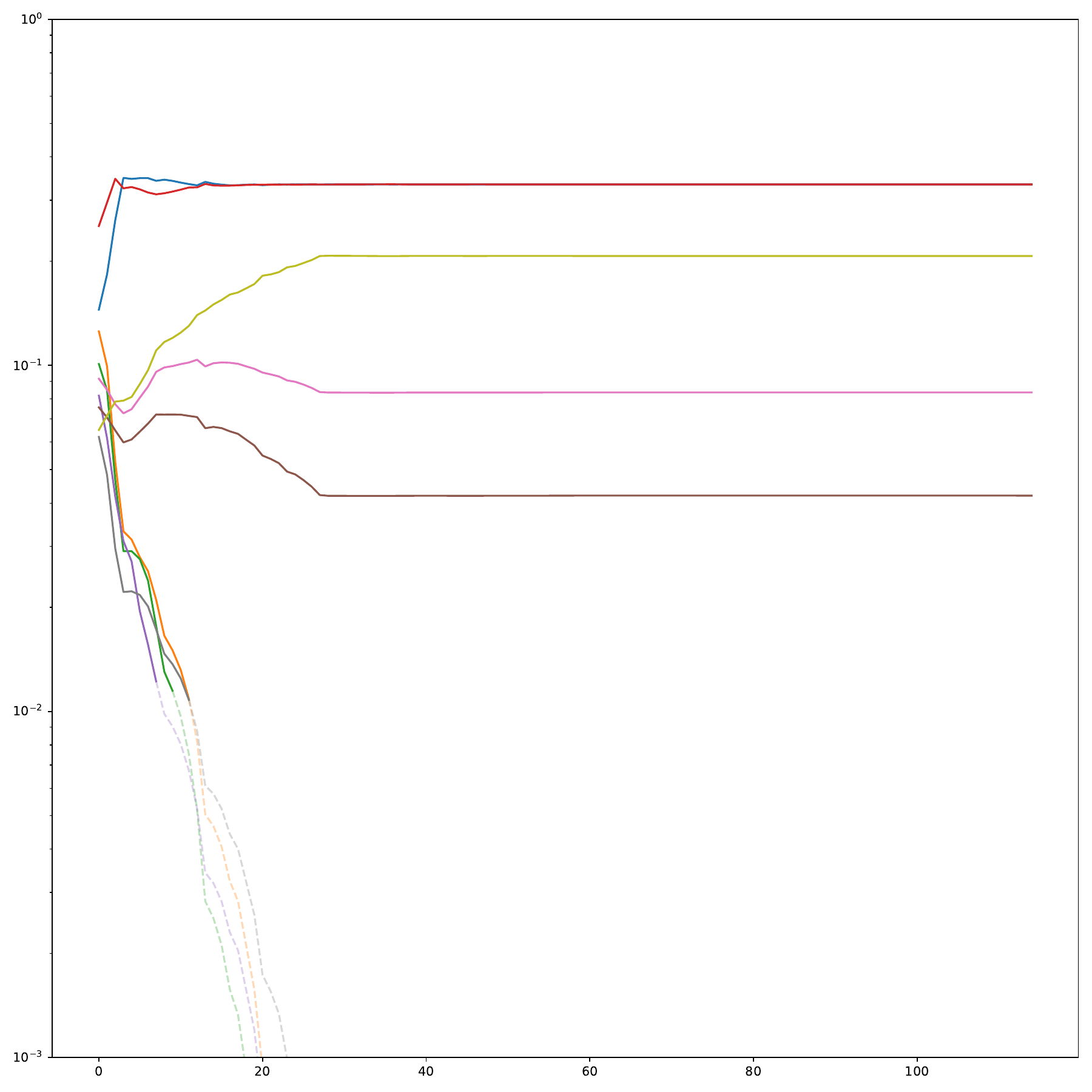}
    \caption{The voting problem. On the left we display the evolution of $\nu_m$ through the iterations, the starting point is a diamond and the ending point is a circle. On the right, we display the value of $p_m$, the probability that the message is sent. When the said probability gets under $1\%$, the curves are dashed and the circle of the corresponding ending $\nu$ is not displayed. The number of messages is respectively 3,6,9 from the top to the bottom and the optimal number of messages is $3$. On the middle, the three extra messages are not used whereas in the bottom, the algorithm ends up with more messages than needed,  the gold, pink and brown messages are the same. The algorithm always end with the correct solution (dirac masses on each middle of the faces).}
    \label{fig:3body problem}
\end{figure}

The optimal solution is obtained for $M=3$ messages. For this problem, Figure~\ref{fig:3body problem} illustrates the evolution of the positions of $\nu_m$ and the corresponding probabilities $p_m$ for different numbers of messages. The algorithm may terminate with an excessive number of messages, either because some messages have a very low probability of occurring (see Figure~\ref{fig:3body problem}, middle or bottom) or because some messages are duplicated (see Figure~\ref{fig:3body problem}, bottom). In \textit{BASIL}, we employ post-processing techniques to decrease the number of messages, ultimately producing a solution that uses the minimal number of messages.

\subsubsection{Numerical convergence of the minimizer with respect to $\varepsilon$}\label{sec:varepsilonsmoothex}

\begin{figure}
    \centering
    \includegraphics[width=0.75\textwidth]{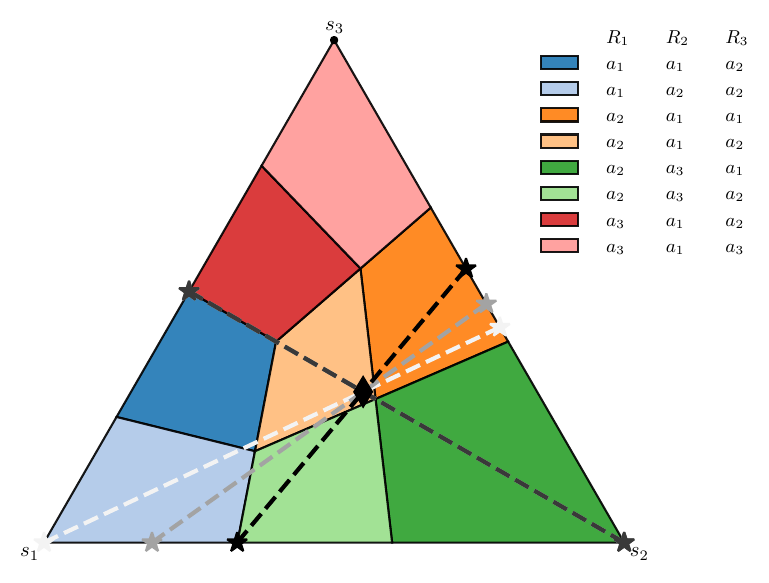}
    \caption{Numerical convergence of the minimizer with respect to $\varepsilon$. For different $\varepsilon$, we compute the optimal signals which are represented as pairs of stars connected by a dashed line (which passes through the prior, represented by a black diamond). Their color depends on the value of $\varepsilon$: the larger the $\varepsilon$, the darker the color.}\label{fig:effectvarepsilon}
\end{figure}

This example, albeit a bit artificial, is chosen to illustrate the convergence of the optimal signals when $\varepsilon$ goes to $0$. There are three states $\mathcal{S} = \{ s_1,s_2,s_3\}$ and three receivers $\{ R_1,R_2,R_3\}$ who need to choose among three actions $\mathcal{A} = \{ a_1,a_2,a_3\}$. The utilities of the receivers are
\begin{equation*}
    u_{r_1} =\begin{pmatrix} 1 & 0 & 0 \\ 0 & 2 & 0\\ 0 & 0 & 1\end{pmatrix},\quad u_{r_2} =\begin{pmatrix} 0 & 1 & 0 \\ 0 & 0 & 2\\ 3 & 0 & 0\end{pmatrix}\quad\text{and}\quad u_{r_3} =\begin{pmatrix} 0 & 3 & 0 \\ 2 & 0 & 0\\ 0 & 0 & 1\end{pmatrix}.
\end{equation*}
The utility of the sender is, for any $s\in\mathcal{S}$,
\begin{align*}
    &v(s,a_1,a_1,a_2) = 2,\; v(s,a_1,a_2,a_2) = 4,\; v(s,a_2,a_1,a_1) = 2,\;v(s,a_2,a_1,a_2) = 1,\\
    & v(s,a_2,a_3,a_1) = 1,\; v(s,a_3,a_1,a_2) = 2\;\text{and}\; \; v(s,a_3,a_1,a_3) = 4.
\end{align*}
The prior is $\mu = (3/10,4/10,3/10)$ and the divergence is the Kullback-Leibler divergence with $\lambda = (1/3,1/3,1/3)$. We compute numerically the optimal signals with two messages for $\varepsilon \in \{10,1,10^{-1}, 10^{-2},10^{-3},10^{-4}\}$ and depict them in Figure~\ref{fig:effectvarepsilon}. As predicted by Theorem \ref{theo:conv:regul}, we observe that the signals vary considerably for large values of $\varepsilon$ before eventually converging.

\subsubsection{The hunter problem}\label{sec:hunterex}

We now consider a scenario in which both the state and action spaces are relatively large. The states represent the locations of two deer within a three-mile square territory. This territory is divided into nine unit-square cells, each identified by coordinates $(i,j)\in\{1,2,3\}^2 = \mathcal{P}$. A state is therefore defined as a pair $s = (s_1,s_2)\in\mathcal{P}^2$, where $s_1$ and $s_2$ indicate the positions of the first and second deer, respectively. 

The receiver in this setting is a hunter who must choose a cell in which to hunt each day, so that his action set is $\mathcal{A} = \mathcal{P}$. The hunter owns a cabin located at some position $p\in\mathcal{P}$, and he prefers to hunt near it. His payoff from choosing action $a$ in state $s$ is
\begin{equation*}
    u(s,a) = h(s,a) - d(p,a)/8,
\end{equation*}
where the function $h$ captures the hunting outcome: it equals $0$ if the chosen cell contains no deer, $1$ if it contains exactly one deer, and $2$ if both deer are present.
The term $d(p,a) = |p_1-a_1| + |p_2-a_2|$ measures the distance between the cabin and the hunting location, thus penalizing choices far from the cabin. The deer's exact positions are not known, but it is observed that they tend to stay close to each other and are attracted to areas with food. The availability of food in each cell is described by a function $f:\mathcal{P}\mapsto \{0,1,2\}$. Combining these behavioral tendencies yields a prior distribution over states:
\begin{equation*}
    \mu(s) = \frac{\nu(s)\kappa(s)}Z,
\end{equation*}
where $Z$ is a normalization constant, $\kappa(s) = 2^{\mathbbm{1}_{\{s_1=s_2\}}}$ reflects the tendency of the deer to remain together and $\nu(s) = 2^{f(s_1)+f(s_2)}$ expresses the influence of food availability.

The second agent, the sender, is a forest ranger who aims to protect a designated conservation zone. This zone consists of three adjacent cells within the territory, formally represented as the set $\mathcal{J} = \{p_k\}_{1\leq k\leq 3}\subset\mathcal{P}$. Each day, the ranger patrols the entirety of the territory and communicates a signal to the hunter in order keep him out of the sanctuary. When the hunter chooses the action $a$, he's utility is equal to $0$ if $a\notin \mathcal{J}$ and $-1$ otherwise. That is, for any state $s$,
\begin{equation*}
    v(s,a) = -\mathbbm{1}_{\mathcal{J}}(a).
\end{equation*}

\begin{figure*}[h!]
    \centering
    \includegraphics[width = 1\textwidth]{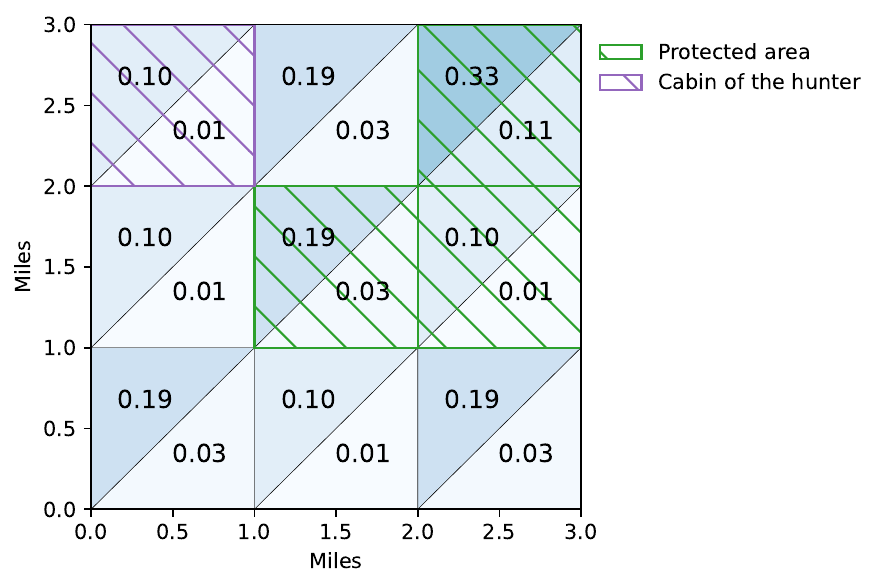}
    \caption{Prior for the hunter problem. The number in the upper-half (resp. the lower-half) of each square is the probability (rounded up) of the presence of a single deer (resp. two deers). Special attention must be given to the upper right cell which is a protected area, close to the cabin and with a high probability of deers.}\label{fig:Hunterprior}
\end{figure*}

With this model, we now turn to numerical simulations with $\varepsilon = 10^{-4}$. For this, we choose $f(q)$ to be equal to one for $q\in\{(1,1),(3,1),(2,2),(2,3)\}$ and two for $q = (3,3)$. The protected area is set as $\mathcal{J} = \{(2,2),(3,2),(3,3)\}$ and the hunter's cabin is located at $p = (1,3)$. The computed prior can be seen in Figure \ref{fig:Hunterprior}.

\begin{figure*}[h!]
    \centering
    \begin{subfigure}[t]{0.5\textwidth}
        \centering
        \includegraphics[width=1\textwidth]{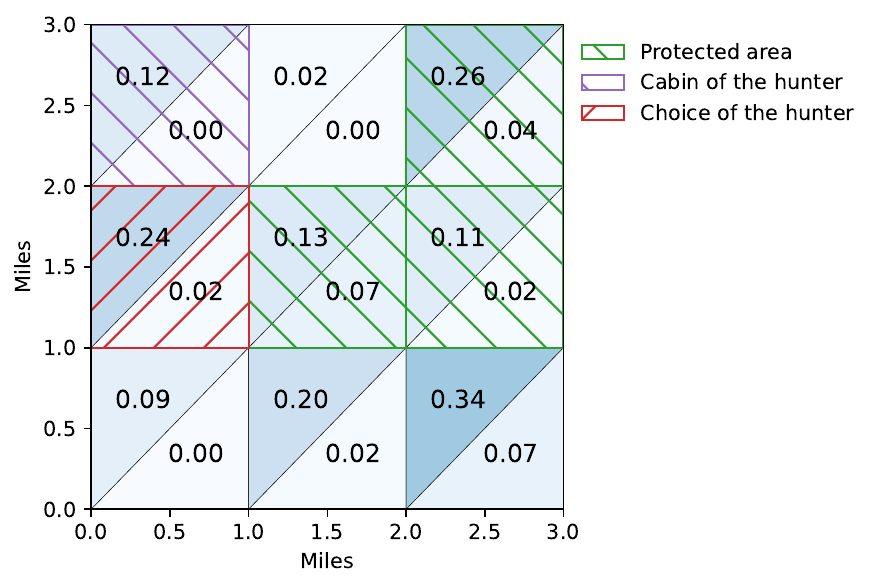}
        \caption{Posterior with first message}
    \end{subfigure}%
    ~ 
    \begin{subfigure}[t]{0.5\textwidth}
        \centering
        \includegraphics[width=1\textwidth]{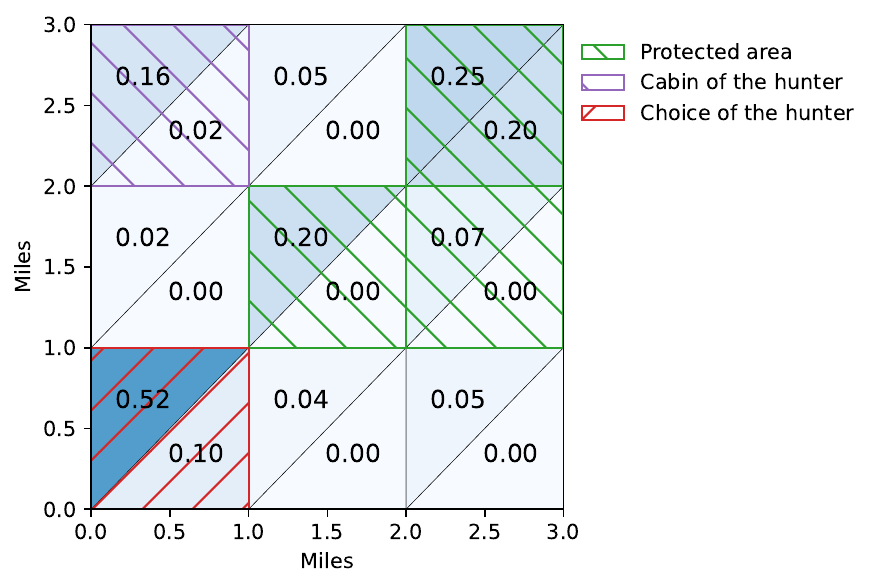}
        \caption{Posterior with second message}
    \end{subfigure}
    \\
    \centering
    \begin{subfigure}[t]{0.5\textwidth}
        \centering
        \includegraphics[width=1\textwidth]{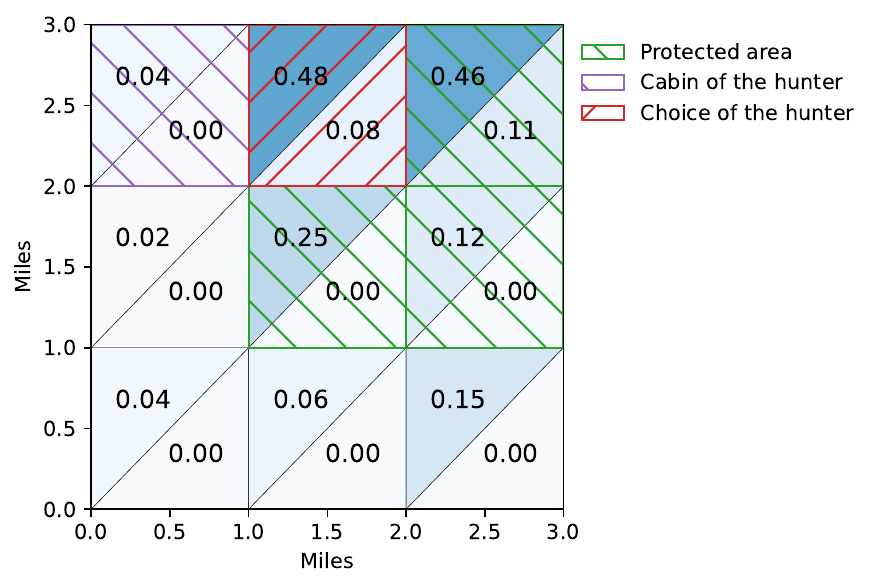}
        \caption{Posterior with third message}
    \end{subfigure}%
    \caption{Depiction of an optimal signal for the hunter problem. The probability (rounded up) of receiving each message is respectively $0.366$ for the first, $0.278$ for the second an $0.356$ for the third. Note that, for each message, the hunter is sent to an unprotected area.}\label{fig:messageshunter}
\end{figure*}

We expect the optimal number of messages to be lower than $|\mathcal{A}| = 9$. It turns out that we can find an optimal signal with $3$ messages as depicted in Figure \ref{fig:messageshunter} where the expected utility of the ranger is (almost) $0$. 

\begin{figure*}[h!]
    \centering
    \includegraphics[width = 1\textwidth]{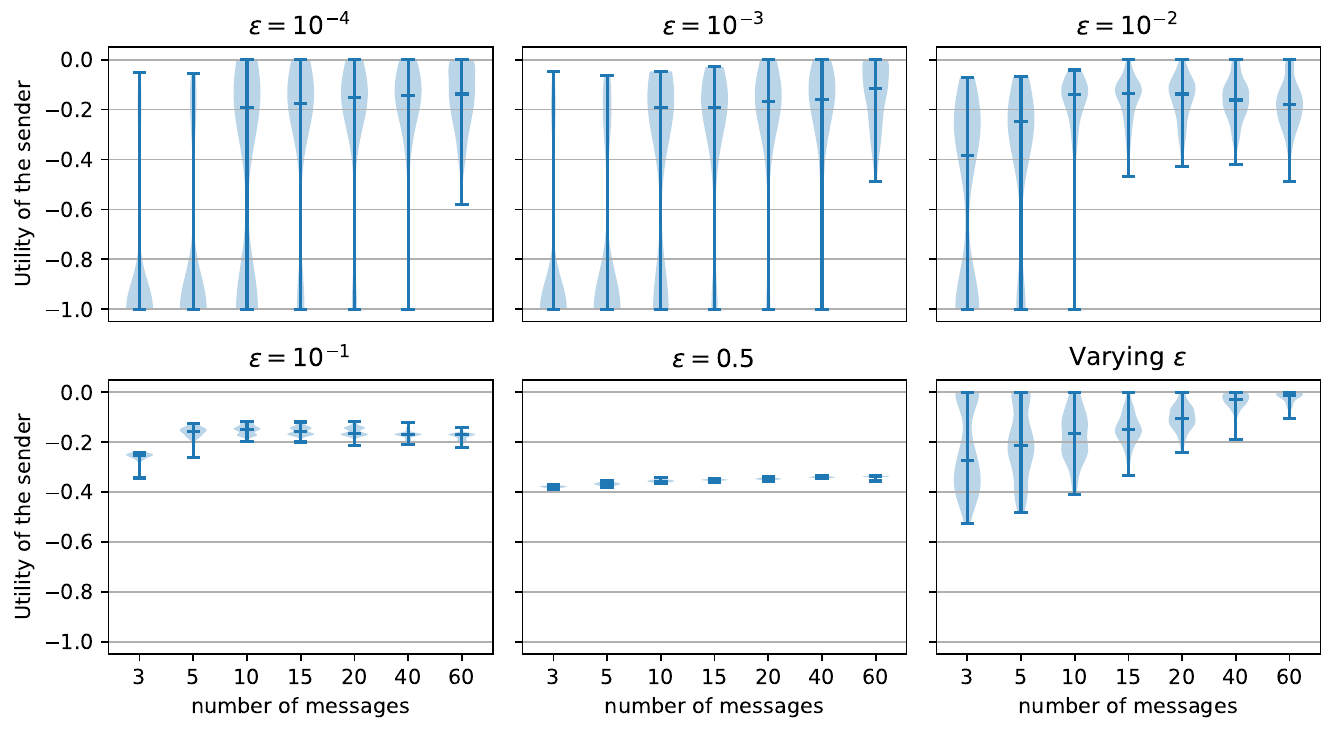}
    \caption{Efficiency of the algorithm for random initial data depending on the values of $\varepsilon$ for the hunter problem.}\label{fig:epsilonhunter}
\end{figure*}

However, because of the randomization of the initial conditions, the algorithm does not consistently converge to the optimal signal. Figure \ref{fig:epsilonhunter} reports the statistics, based on one hundred realizations, of the expected utility of the sender for randomized initial data under varying values of $\varepsilon$. As the figure illustrates, the regularization process tends to obscure information from the optimization problem, so that the optimal value is attained only when $\varepsilon$ is sufficiently small. However, for small values of $\varepsilon$, the algorithm tends to yield (very) bad signals when the number of messages in small. We observe that increasing the number of messages improves the convergence of the algorithm toward signals that more closely approximate the optimal one. Our interpretation is that artificially enlarging the message space strengthens the algorithm's exploratory and selective capacity (has also seen in Section \ref{sec:votingproblem}), thereby facilitating convergence to the optimal set of messages. Based on this observation, we iteratively solve the optimization problem for decreasing values of $\varepsilon$ of the form $10^{-k}$, where $k$ is linearly spaced between $0$ and $4$, yielding a total of ten distinct values. This gives a more robust and effective approach whose results can be seen in the last statistic, titled "Varying $\varepsilon$". This approach is natively implemented in the \textit{BASIL} library.

%\subsubsection{Numerical examples TODO}
%Examples du vote, influence du nombre de messages, cela aide à la convergence,
%influence de $\varepsilon$, si c'est trop petit, les gradients s'annulent et ce n'est pas utile, il vaut mieux avoir de $\varepsilon$ grand.
%
%Capacité exploratoires du problème quand  tu pars de  points qui engendrent la même action, points proches mais points éloignés posssibles ?
%
%Example avec plusieurs types votes avec (2000) . Incertitude sur le receiver.

\section{Conclusion}

This article presents a new approach to solving the Bayesian persuasion problem based on regularization methods. Our method has the advantage of ensuring the receiver's solution is unique (and explicit in certain cases), which makes it possible to employ first- and second-order optimization methods. We prove that the solution of the regularized problem converges to that of the original problem as $\varepsilon$ tends to $0$. In addition, we provide a version of the revelation principle that allows one to determine the optimal number of messages for a given problem. Through various numerical examples, we examine the strengths and limitations of our method and justify our numerical choices. For these experiments, we developed a Python library, \textit{BASIL}, which is publicly available and ensures the reproducibility of the reported results.

\section{Proofs}

\subsection{Proof of Lemma \ref{lem:compactau}}

We know that the space $\mathsf{P}(\mathcal{S})$ endowed with the Prokhorov metric is a compact space so that the space $\mathsf{P}(\mathsf{P}(\mathcal{S}))$ (also endowed with the Prokhorov metric) is a compact space. Let $\mu\in\mathsf{P}(\mathcal{S})$ and $(\tau_n)_{n\in\mathbb{N}}$ be a sequence in $\mathsf{T}_{\mu}$. We know that, up to a subsequence,  $(\tau_n)_{n\in\mathbb{N}}$ converges to an element $\tau\in \mathsf{P}(\mathsf{P}(\mathcal{S}))$. Denote
\begin{equation*}
    \eta = \int_{\mathsf{P}(\mathcal{S})} \nu d\tau(\nu). 
\end{equation*}
We wish to prove that $\eta = \mu$. For any $s\in\mathcal{S}$, we consider the continuous bounded function $f_s: \mathsf{P}(\mathcal{S})\mapsto [0,1]$ given by
\begin{equation*}
    f_s(\nu) = \nu(s).
\end{equation*}
Thus, we know that, for any $s\in\mathcal{S}$, up to a subsequence,
\begin{equation*}
    \mu(s) = \int_{\mathsf{P}(\mathcal{S})}f_s(\nu)d\tau^n(\nu) \underset{n\to\infty}\longrightarrow \int_{\mathsf{P}(\mathcal{S})}f_s(\nu)d\tau(\nu) = \nu(s),
\end{equation*}
which proves that $\tau\in\mathsf{T}_{\mu}$ and, thus, $\mathsf{T}_{\mu}$ is compact. 

Let $M\in\mathbb{N}^*$, $\mu\in\mathsf{P}(\mathcal{S})$ and $(\tau_n)_{n\in\mathbb{N}}$ be a sequence in $\mathsf{T}_{M,\mu}$. We have that
\begin{equation*}
    \tau_n = \sum_{m = 1}^M p_{m,n}\delta_{\nu_{m,n}}.
\end{equation*}
Since, for any $n\in\mathbb{N}$, $(p_{m,n},\nu_{m,n})_{1\leq m\leq M} \in ([0,1]\times\mathsf{P}(\mathcal{S}))^M$, where $([0,1]\times\mathsf{P}(\mathcal{S}))^M$ is a compact space, we know that there exists $(p_{m},\nu_{m})_{1\leq m\leq M} \in ([0,1]\times\mathsf{P}(\mathcal{S}))^M$ such that, up to a subsequence,
\begin{equation*}
    \tau_n \underset{n\to\infty}{\longrightarrow} \sum_{m = 1}^M p_m\delta_{\nu_{m}} \in\mathsf{T}_{M,\mu},
\end{equation*}
which proves the desired result.

\subsection{Proof of Lemma \ref{lem:Wusc}}

We denote $\mathrm{co}(W)$ the convex hull of $\mathrm{Gr}(W) = \{(\nu,W(\nu))\;\text{for any }\nu\in\mathsf{P}(\mathcal{S})\}$. That is, for any $\nu\in\mathsf{P}(\mathcal{S})$ and $(\nu,v)\in \mathrm{co}(W)$, there exists $\tau\in\mathsf{T}_{\nu}$ such that
\begin{equation*}
    v = \int_{\mathsf{P}(\mathcal{S})}W(\nu)d\tau(\nu).
\end{equation*}
In other words, $\mathrm{co}(W)$ is also given by
\begin{equation*}
    \mathrm{co}(W) = \left\{\left(\mu,\int_{\mathsf{P}(\mathcal{S})}W(\nu)d\tau(\nu)\right)\text{ for any }\mu\in\mathsf{P}(\mathcal{S})\text{ and }\tau\in\mathsf{T}_{\mu}\right\}.
\end{equation*}
We define, for any $\nu\in\mathsf{P}(\mathcal{S})$,
\begin{equation*}
    \tilde{W}(\nu) = \sup\{v\text{ such that }(\nu,v)\in\mathrm{co}(W)\},
\end{equation*}
so that, if $W$ is upper semicontinuous, we have $\tilde{W} = \hat{W}$.
The proof is decomposed in several points.

\subsubsection{$W$ is upper semicontinuous}

Let $(\nu_{n})_{n\in\mathbb{N}}\subset\mathsf{P}(\mathcal{S})$ be a sequence that converges to $\nu$ as $n\to\infty$. For any $t\in\mathcal{T}$, since $(a^{\star}_{t}(\nu_n))_{n\in\mathbb{N}}\subset\mathcal{A}$, we know that it converges, up to a subsequence, to a $a^*_t\in\mathcal{A}$.
For any $\ell\in\{1,\ldots,L\}$ and $t\in\mathcal{T}$, we have, for any $a_{\ell}\in\mathcal{A}_{\ell}$,
\begin{equation*}
    \sum_{s\in\mathcal S} u_{t_\ell}(s,a_{\ell})\nu_n(s) \leq \sum_{s\in\mathcal S} u_{t_\ell}(s,(a^{\star}_{t}(\nu_n))_{\ell})\nu_n(s).
\end{equation*}
Thus, by passing to the limit, up to a subsequence, in the previous inequality yields that $a^*_t\in\mathcal{A}_{t}(\nu)$. It follows that, up to a subsequence,
\begin{align*}
    W(\nu_n) &= \sum_{(s,t)\in\mathcal{S}\times\mathcal{T}}v(s,a_t^\star(\nu_n)) \nu_n(s) \eta(t) 
    \\&\underset{n\to\infty}{\longrightarrow} \sum_{(s,t)\in\mathcal{S}\times\mathcal{T}}v(s,a^*_t) \nu(s) \eta(t) 
    \\ &\hspace{3em}\leq \sum_{(s,t)\in\mathcal{S}\times\mathcal{T}}v(s,a_t^\star(\nu)) \nu(s) \eta(t) = W(\nu) ,
\end{align*}
so that $\mathrm{limsup}_{n\to\infty} W(\nu_n)\leq W(\nu)$, which is the desired result.

\subsubsection{The function $\tilde{W}$ is concave}

Indeed, for any $\nu_1,\nu_2\in\mathsf{P}(\mathcal{S})$ and $r\in[0,1]$, we notice that
\begin{equation*}
    \left(r\nu_1 + (1-r)\nu_2,rW(\nu_1) + (1-r)W(\nu_2)\right)\in\mathrm{co}(W).
\end{equation*}
Let $(W_1^n)_{n\in\mathbb{N}}$ and $(W_2^n)_{n\in\mathbb{N}}$ such that $(\nu_1,W_1^n),(\nu_2,W_2^n)\in\mathrm{co}(W)$ and
\begin{equation*}
    W_1^n\underset{n\to\infty}{\longrightarrow}\tilde{W}(\nu_1)\quad\text{and}\quad W_2^n\underset{n\to\infty}{\longrightarrow}\tilde{W}(\nu_2).
\end{equation*}
We have, for any $n\in\mathbb{N}$,
\begin{equation*}
    \left(r\nu_1 + (1-r)\nu_2,rW_1^n + (1-r)W_2^n\right)\in\mathrm{co}(W),
\end{equation*}
so that
\begin{equation*}
    rW_1^n + (1-r)W_2^n\leq \tilde{W}(r\nu_1 + (1-r)\nu_2),
\end{equation*}
which yields, by passing to the limit $n\to\infty$,
\begin{equation*}
    r\tilde{W}(\nu_1) + (1-r)\tilde{W}(\nu_2)\leq \tilde{W}(r\nu_1 + (1-r)\nu_2).
\end{equation*}
Thus, $\tilde{W}$ is concave.

\subsubsection{The function $\tilde{W}$ is the smallest concave function greater that $W$}
Let $V:\mathsf{P}(\mathcal{S})\mapsto \mathbb{R}$ be a concave function. On one hand, for any $(\nu,v)\in\mathrm{co}(V)$, there exists $\tau\in\mathsf{T}_{\nu}$ such that
\begin{equation*}
    v = \int_{\mathsf{P}(\mathcal{S})} V(\nu)d\tau(\nu) \leq V(\nu),
\end{equation*}
where we used the fact that $V$ is concave. It follows that $\tilde{V}(\nu)\leq V(\nu)$. On the other hand, since $\mathrm{Gr}(V)\subset\mathrm{co}(V)$, we have that $V(\nu)\leq \tilde{V}(\nu)$. We conclude that $\tilde{V} = V$ if $V$ is concave.

Now, let $H:\mathsf{P}(\mathcal{S})\mapsto \mathbb{R}$ be such that $W\leq H$. For any $(\nu,w)\in\mathrm{co}(W)$, there exists $\tau\in\mathsf{T}_{\nu}$ such that
\begin{equation*}
   w =  \int_{\mathsf{P}(\mathcal{S})} W(\nu)d\tau(\nu)\leq \int_{\mathsf{P}(\mathcal{S})} H(\nu)d\tau(\nu),
\end{equation*}
and, in particular, we have $\tilde{H}(\nu)\geq w$ which yields $\tilde{W}(\nu)\leq \tilde{H}(\nu)$.

It follows that, for any $V:\mathsf{P}(\mathcal{S})\mapsto \mathbb{R}$ which is concave and such that $W\leq V$, we have $\tilde{W}\leq V$. Thus, $\tilde{W}$ is the smallest concave function greater than $W$.

\subsection{Proof of Lemma \ref{lem:decompotau}}

We proceed by contradiction and assume that there exists $\tilde{\tau_1}\in\mathsf{T}_{\mu_1}$ such that
\begin{equation*}
    \int_{\mathsf{P}(\mathcal{S})}W(\nu)d\tau_1(\nu)< \int_{\mathsf{P}(\mathcal{S})}W(\nu)d\tilde{\tau}_1(\nu).
\end{equation*}
Then, by setting $\tilde{\tau} = r\tilde{\tau_1} + (1-r)\tau_2$, we observe that $\tilde{\tau}\in\mathsf{T}_{\mu}$ and
\begin{align*}
    &\int_{\mathsf{P}(\mathcal{S})}W(\nu)d\tilde{\tau}(\nu) = r\int_{\mathsf{P}(\mathcal{S})}W(\nu)d\tilde{\tau}_1(\nu) + (1-r)\int_{\mathsf{P}(\mathcal{S})}W(\nu)d\tau_2(\nu) 
    \\ &\hspace{3em}> r\int_{\mathsf{P}(\mathcal{S})}W(\nu)d\tau_1(\nu)  + (1-r)\int_{\mathsf{P}(\mathcal{S})}W(\nu)d\tau_2(\nu)  = \int_{\mathsf{P}(\mathcal{S})}W(\nu)d\tau(\nu) ,
\end{align*}
which contradicts the fact that $\tau$ is a solution of \eqref{eq:senderprob:tausupp}. The same arguments hold for $\tau_2$.

\subsection{Proof of Corollary \ref{cor:concavWtau}}

By Remark \ref{rmk:representerthm}, we know that $\tau$ can be decomposed as
\begin{equation*}
    \tau = \sum_{m = 1}^M p_m\delta_{\nu_m},
\end{equation*}
for some $M\leq |\mathcal{S}|$, $(p_m)_{1\leq m\leq M}\in\Delta_M$ and $(\nu_m)_{1\leq m\leq M}\subset\mathsf{P}(\mathcal{S})$ verifying $\sum_{m =1}^M p_m\nu_m = \mu$. We observe that we can decompose $\tau = r\tau_1 + (1-r)\tau_2$ with
\begin{equation*}
    \tau_1 = \delta_{\nu_1},\quad \tau_2 = \sum_{m = 2}^M p_m\delta_{\nu_m}\quad\text{and}\quad r= p_1.
\end{equation*}
With have $\tau_1\in\mathsf{T}_{\nu_1}$ and, by Lemma \ref{lem:decompotau}, $\tau_1$ is a solution of \eqref{eq:senderprob:tau}. Thus, we have
\begin{equation*}
    W(\nu_1) = \int_{\mathsf{P}(\mathcal{S})} W(\nu)d\tau_1(\nu) = \hat{W}(\nu_1).
\end{equation*}
By iterating this argument, we obtain that $W(\nu_m)= \hat{W}(\nu_m)$, for any $1\leq m\leq M$.

\subsection{Proof of Theorem \ref{thm:tausupp}}
\subsubsection{Upper bound on the number of messages}

We now discuss the number of possible actions. We start with the first receiver $\ell$ and for each type $t_{\ell}$, we choose its action in $\mathcal A_\ell$, there are  $A_\ell^{T_\ell}$ choices and iterating over the receivers, the set of available actions is $\mathcal A^\sharp=\bigotimes_{\ell \in \mathcal L}\bigotimes_{t_\ell \in \mathcal T_\ell} \mathcal A_\ell$ that can be rewritten as:
\[\mathcal A^\sharp =\{a^\sharp=((a^{\sharp}_{t_\ell})_{t_\ell \in \mathcal T_\ell})_{\ell\in \mathcal L}
\text{ such that } a_{t_\ell} \in {\mathcal A_\ell } \quad\forall \ell \in \mathcal L ,t_\ell \in \mathcal T_\ell \}\]
This set is of cardinal $\prod_{\ell \in \mathcal L}A_\ell^{T_\ell}$ and can be rewritten as .

We tensorize the action by types, that is we consider a new set of actions $\mathcal A^\sharp=\bigotimes_{t\in \mathcal T} \mathcal A $ of cardinal $A^T$ which is defined as 

For any $a^\mathcal T \in \mathcal A^\mathcal T$, we introduce the set $\mathcal C(a^\mathcal T)$ which is the set of priors $\nu$ such that, for each type $t$ the action $a^\mathcal T_{t}$ is admissible for the receiver $\ell$ of type $t_\ell$, that is 
\[\mathcal C(a^\mathcal T)=\{\nu \in \mathrom P(\mathcal S) \text{ such that, for each } t\in \mathcal T, \ell \in \mathcal L, \text{ then } (a^\mathcal T_{t})_\ell \in \mathcal A^\star_{t_\ell}(\nu)\}\]
The set $\mathcal C(a^\mathcal T)$ is convex, indeed $\nu$ belongs to $\mathcal C(a^\mathcal T)$ if and only if we have for every $t\in \mathcal T, \ell \in \mathcal L$ and forall $a\in \mathcal A_\ell$,
\[\sum_s u_{t_\ell}(s,a)\nu(s) \le \sum_s u_{t_\ell}(s,(a^\sharp_{t})_\ell)\nu(s).\]
If the above inequality is true for $\nu_1$ and $\nu_2$, it is surely true for $r\nu_1+(1-r)\nu_2$ for any $r\in [0,1]$, and this proves that $\mathcal C(a^\mathcal T)$ is convex.

Let $\tau\in\mathrom{T}_\mu$ be a solution of \eqref{eq:senderprob:tau} and for each $a^\sharp \in \mathcal A^\sharp$ denote $\nu(a^\sharp)$, the average of the prior $\nu$ over the set of prior that admits $a^\sharp$ as action policies.
\[p(a^\sharp) =\int_{a^\star(\nu)=a^\sharp}d\tau(\nu)
\text{ and } \nu(a^\sharp)=\frac{\int_{a^\star(\nu)=a^\sharp}\nu d\tau(\nu)}{p(a^\sharp)}.
\]
Note that, thanks to the discussion in preamble of this section, we restrict our analysis to the set $A^\sharp \subset A^\mathcal T$ and $\mathrom{P}(S)$ is partitionned into the different sets $\{\nu, a^\star(\nu)=a^\sharp\}_{a^\sharp \in \mathcal A^\sharp}$. 
Introducing $\tilde \tau=\sum_{a^\sharp \in \mathcal A^\sharp} p(a^\sharp)\delta_{\nu(a^\sharp)}$, we have
\begin{eqnarray*}
\int_{\mathrom P(S)} \nu d\tilde \tau(\nu)
=\sum_{a^\sharp \in \mathcal A^\sharp}
p(a^\sharp)\nu(a^\sharp)
=\sum_{a^\sharp \in \mathcal A^\sharp}
\int_{a^\star(\nu)=a^\sharp}\nu d\tau(\nu)
=\int_{\mathrom P(S)}\nu d\tau(\nu)
=\mu
\end{eqnarray*}
Similarly, we prove that $\sum_{a^\sharp}p(a^\sharp)=1$, so that $\tilde \tau$ belongs to $\mathrom T_\mu$ and is admissible in Problem~\eqref{eq:senderprob:tau}.
Because $a^\star(\nu)$ is chosen amongst the admissible actions, then 
\[ a^\star(\nu)=a^\sharp \Rightarrow \nu \in \mathcal C(a^\sharp).\]
By definition, $\nu(a^\sharp)$ is an average of priors $\nu$ which all belong to the convex set $\mathcal C(a^\sharp)$, hence $\nu(a^\sharp)\in \mathcal C(a^\sharp)$. Hence 
\[W(\nu(a^\sharp))=\sum_{s,t} v(s,a^\star_t(\nu(a^\sharp)))\nu(a^\sharp) \eta(t) \ge 
\sum_{s,t} v(s,a^\sharp_t)\nu(a^\sharp) \eta(t)
\]

It follows that
\begin{eqnarray*}
\int_{\mathcal P(\mathcal S)} W(\nu)d\tilde \tau(\nu)
&=&
\sum_{a^\sharp} p(a^\sharp)  W(\nu(a^\sharp))
\ge
\sum_{a^\sharp} \sum_{s,t} v(s,a^\sharp)p(a^\sharp)\nu(a^\sharp) \eta(t)\\
&\ge&\sum_{a^\sharp} \sum_{s,t} v(s,a^\sharp)\left(\int_{a^\star(\nu)=a^\sharp}\nu d\tau(\nu) \eta(t)\right)\\
&\ge& \int_{\mathrom P(\mathcal S)} \sum_{s,t} v(s,a^\star(\nu))\nu d\tau(\nu) \eta(t)=\int_{\mathrom P(\mathcal S)} W(\nu)d \tau(\nu)
\end{eqnarray*}
Because $d\tau$ is optimal, so is $d\tilde \tau$ and every inequality becomes an equality and $a^\sharp$ is not only admissible for $\nu(a^\sharp)$ but can be defined as the action policy taken by the receivers. That is, we can suppose that $a^\sharp=a^\star_t(\nu(a^\sharp))$. Hence, as claimed, each Dirac measure that compose the optimal $d\tilde \tau$ is associated to a different action policy $a^\sharp$. Finally there is at most $\prod_\ell A_\ell^{\mathcal T_\ell}$ of them.
\subsubsection{Lower bound on the number of messages}
We are given $\mathcal A$ and we suppose that we have only one type, that is $T=1$. We want to design utility functions $u$ and $v$ on a precise space state $\mathcal S$ so that there is no optimal solution to \eqref{eq:senderprob:tau} with support of cardinality lower than or equal to $A$. For that purpose, we suppose that $\mathcal S=\mathcal A$.
For each $s=a \in \mathcal S$ we define the utility of the receivers and the one of the sender to be
\[u_{t_\ell}(s,a_\ell)=\begin{cases} 1 &\text{ if } a_\ell=s_\ell \\ 0 & \text{ if } a_\ell \in \mathcal A_\ell \setminus \{s_\ell\}\ 
\end{cases}
\quad\text{and}\quad v(s,a)=\begin{cases} 1 &  s=a \\ 
0 & \text{ if not} \end{cases}.\]
Clearly, for each $s\in \mathcal S$, if $\nu_s=\delta_s$, then there is only one admissible action which is $s$ for which the gain of the sender is $W(\nu_s)=1$. Take any $\mu=\sum_{s\in \mathcal S} \mu(s) \nu_s$ which is not located on the vertices of $\mathrom P(\mathcal S)$, that is its support is of cardinal at least $2$, or equivalently $\mu(s)<1$ for every $s$, recall that, for this particular $\mu$, the receiver chose a certain action $a^\star(\mu)$ and we have, by construction of $v$,
\[W(\mu)=\sum_s v(s,a^\star(\mu))\mu(s)= \mu(a^\star(\mu))\]
Consider now $\tau = \sum_{s} \mu(s)\delta_{\nu_s}$, where $\nu_s=\delta_s$. It is easy to check that $\tau\in \mathrom{T}_{\mu}$ and then
\[\hat W(\mu)\ge \int_{\mathrom P(\mathcal S)}W(\nu)d\tau(\nu)=\sum_s \mu(s) W(\nu_s).\]
As we said before the sender gains $1$ for $\nu_s$, so that we obtain
\[\hat W(\mu)\ge \sum_s \mu(s) > \mu(a^\star(\mu))=W(\mu)\]
So that for every measure $\mu$ which is not of the form $\mu=\nu_s$, we have $\hat W(\mu)>W(\mu)$.
Now take any $\mu$ such that $\mu(s) >0$ for every $s$ and take $\tau$ optimal for the problem of finding $\hat W(\mu)$. If $\tau$ is of finite support with $M$ Dirac masses, then each Dirac mass must be supported on a $\nu_s$ such that $\nu_s=\delta_s$ for some $s\in \mathcal S$. And then $\tau$ must be written as $\tau =\sum_s \tau_s \delta_{\nu_s}$ and the condition $\int \nu \tau(\nu)\mu$ imposes $\tau_s=\mu(s)\ne 0$ for every $s$ and hence $\tau$ has support of cardinality $S=A$.

\subsection{Proof of Lemma \ref{lem:limsupvaluefun}}

We denote
\begin{equation*}
    \tau^{\varepsilon} = \sum_{m = 1}^M p_m^\varepsilon\delta_{\nu_m^\varepsilon},
\end{equation*}
where $(p_m^\varepsilon)_{1\leq m\leq M}\in\Delta_M$ and $(\nu_m^\varepsilon)_{1\leq m\leq M}\subset\mathsf{P}(\mathcal{S})$ verifying $\sum_{m =1}^M p_m^\varepsilon\nu_m^\varepsilon = \mu$. We also denote, for any $\ell\in\{1,\ldots,L\}$ and any $t\in\mathcal{T}_{\ell}$,
\begin{equation*}
    \theta^{\star,\varepsilon}_{m} = (\theta^{\star,\varepsilon}_{\ell,t_{\ell},\nu_m^\varepsilon})_{\substack{1\leq \ell\leq L\\t\in\mathcal{T}^L}}\in\bigoplus_{\substack{1\leq \ell\leq L\\t\in\mathcal{T}^L}}\mathsf{P}(\mathcal{A}_{\ell}) = \mathfrak{X}.
\end{equation*}
Since $\mathfrak{X}^M$ is a compact space, we know that, up to a subsequence,
\begin{equation*}
    (p_m^{\varepsilon},\nu_m^{\varepsilon},\theta^{\star,\varepsilon}_m)_{1\leq m\leq M}\underset{\varepsilon\to0}{\longrightarrow}(p_m,\nu_m,\theta^{\star}_m)_{1\leq m\leq M},
\end{equation*}
with
\begin{equation*}
    \tau^0 = \sum_{m = 1}^M p_m\delta_{\nu_m}.
\end{equation*}
We denote, for any $1\leq m\leq M$,
\begin{equation*}
     (\theta^{\star}_{\ell,t_{\ell},m})_{\substack{1\leq \ell\leq L\\t\in\mathcal{T}^L}} = \theta^{\star}_m.
\end{equation*}
Since, for any $\ell\in\{1,\ldots,L\}$ and any $t_{\ell}\in\mathcal{T}_{\ell}$, $\theta^{\star,\varepsilon}_{\ell,t_{\ell},\nu_m^\varepsilon}$ is a solution of \eqref{eq:probregreceiver:util}, we have, for every $\theta\in\mathsf{P}(\mathcal{A}_{\ell})$,
\begin{align*}
    &\sum_{(s,a)\in\mathcal{S}\times\mathcal{A}_\ell} u_{t_{\ell}}(s,a)\nu_m^\varepsilon(s)\theta(a) - \varepsilon \vartheta(\theta,\lambda_{\ell})
    \\ &\hspace{5em}\leq\sum_{(s,a)\in\mathcal{S}\times\mathcal{A}_\ell} u_{t_{\ell}}(s,a)\nu_m^\varepsilon(s)\theta^{\star,\varepsilon}_{\ell,t_{\ell},\nu_m^\varepsilon}(a) - \varepsilon \vartheta(\theta^{\star,\varepsilon}_{\ell,t_{\ell},\nu_m^\varepsilon},\lambda_{\ell}) 
    \\ &\hspace{10em}\leq \sum_{(s,a)\in\mathcal{S}\times\mathcal{A}_\ell} u_{t_{\ell}}(s,a)\nu_m^\varepsilon(s)\theta^{\star,\varepsilon}_{\ell,t_{\ell},\nu_m^\varepsilon}(a).
\end{align*}
Thus, by passing to the limit $\varepsilon\to 0$ (of the subsequence), we deduce that
\begin{equation*}
    \max_{\theta\in\mathsf{P}(\mathcal{A}_{\ell})}\sum_{(s,a)\in\mathcal{S}\times\mathcal{A}_\ell} u_{t_{\ell}}(s,a)\nu_m(s)\theta(a) = \sum_{(s,a)\in\mathcal{S}\times\mathcal{A}_\ell} u_{t_{\ell}}(s,a)\nu_m(s)\theta^{\star}_{\ell,t_{\ell},m}(a),
\end{equation*}
so that $\theta^{\star}_{\ell,t_{\ell},m}$ belongs in $\Theta_{\ell,t_\ell} (\nu_m)$. In particular, by denoting $\theta^{\star}_{t,m} = \prod_{\ell = 1}^L \theta^{\star}_{\ell,t_{\ell},m}$, we observe that $\theta^{\star}_{t,m}\in\Theta_{t} (\nu_m)$ and, thus,
\begin{align*}
    \sum_{(s,a)\in\mathcal{S}\times\mathcal{A}}v(s,a)\theta_{t,m}^\star(a)\nu_m(s) &\leq \max_{\theta\in\Theta_{t} (\nu_m)}\sum_{(s,a)\in\mathcal{S}\times\mathcal{A}}v(s,a)\theta(a)\nu_m(s)
    \\ &\leq \sum_{(s,a)\in\mathcal{S}\times\mathcal{A}}v(s,a)\theta^{\star}_{t,\nu_m}(a)\nu_m(s).
\end{align*}
Hence, by denoting, for any $a\in\mathcal{A}^{\mathcal{T}}$,
\begin{equation*}
    \theta^{\star}_{m}(a) = \prod_{t\in\mathcal{T}^L}\theta^{\star}_{t,m}(a_t)\quad\text{and}\quad W^0_m = \sum_{(s,a)\in\mathcal{S}\times\mathcal{A}^{\mathcal{T}}} w(s,c)\nu_m(s)\theta^{\star}_{m}(a),
\end{equation*}
it follows that, up to a subsequence,
\begin{equation*}
    \mathsf{W}^{\varepsilon}(\tau^{\varepsilon}) = \sum_{m = 1}^M W^{\varepsilon}(\nu_m^{\varepsilon})p_m^{\varepsilon}\underset{\varepsilon\to0}{\longrightarrow} \sum_{m = 1}^MW^0_m p_m \leq \sum_{m = 1}^MW(\nu_m) p_m = \mathsf{W}(\tau^0),
\end{equation*}
which yields the desired result.

\subsection{Proof of Lemma \ref{lem:convergtaureg}}

Let $\tau^{\star}$ be a solution of \eqref{eq:senderprobreg:tau}. Theorem \ref{thm:tausupp} yields the existence of a set $\mathcal{R}\subset\mathcal{A}^{\mathcal{T}}$ such that
\begin{equation*}
    \tau^{\star} = \sum_{a\in\mathcal{R}}p_a^{\star}\delta_{\nu_a^{\star}},
\end{equation*}
with $\nu_a^{\star}$ such that $a_t = a^{\star}_t(\nu_a^{\star})$, for any $t\in\mathcal{T}$.
Since, by Assumption \ref{asm:reg}, for any $a\in\mathcal{R}$, there exists $\nu_a$ such that $\mathrm{supp}(\nu_a)\subset\mathrm{supp}(\mu)$ and $\mathcal{A}_t(\nu_a) = \{a_t\}$, for any $t\in\mathcal{T}$. Furthermore, since $\mathcal{S}$ is finite, there exists $r>0$ such that
\begin{equation*}
    \bar{\nu} = \mu + r\left(\mu - \sum_{a\in\mathcal{R}}p_a^{\star}\nu_a\right)\in\mathsf{P}(\mathcal{S}).
\end{equation*}
Let $\alpha = \varepsilon^{1/2}$ as well as, for any $a\in\mathcal{R}$,
\begin{align*}
    \nu^{\varepsilon}_a = \alpha \nu_a + (1-\alpha)\nu_a^{\star},\; \beta = \frac{\alpha}{\alpha+r}\quad\mbox{and}\quad \tau^{\varepsilon} = (1-\beta)\sum_{a\in\mathcal{R}}p_a^{\star}\delta_{\nu_a^{\varepsilon}} + \beta\delta_{\bar{\nu}}.
\end{align*}
We have
\begin{equation*}
    (1-\beta)\sum_{a\in\mathcal{R}}p_a + \beta = 1,
\end{equation*}
so that $\tau^{\varepsilon}\in\mathsf{P}(\mathsf{P}(\mathcal{S}))$ and, furthermore, we observe that
\begin{align*}
    \int_{\mathsf{P}(\mathcal{S})}\nu d\tau^{\varepsilon}(\nu) &= (1-\beta)\sum_{a\in\mathcal{R}}p_a^{\star}\nu_a^{\varepsilon} + \beta\bar{\nu}
    \\ &=(1-\beta)\sum_{a\in\mathcal{R}}p_a^{\star}\nu_a^{\varepsilon} + \beta\left(\mu + r\left(\mu - \sum_{a\in\mathcal{R}}p_a^{\star}\nu_a\right)\right)
    \\ &= \beta(1+r)\mu + \sum_{a\in\mathcal{R}}p_a^{\star}\left((1-\beta)\nu_a^{\varepsilon} - \beta r \nu_a\right)
    \\ &= \beta(1+r)\mu + (1-\beta)(1-\alpha)\sum_{a\in\mathcal{R}}p_a^{\star}\nu_a^{\star}
    \\&=\left(\beta(1+r) + (1-\beta)(1-\alpha) \right)\mu = \mu,
\end{align*}
which yields that $\tau^{\varepsilon}\in\mathsf{T}_{\mu}$. Since $\mathcal{A}_t(\nu_a) = \{a_t\}$, we have, for any $\ell\in\{1,\ldots,L\}$ and $t_{\ell}\in\mathcal{T}_{\ell}$, that $\mathcal{A}_{\ell,t_{\ell}}(\nu_a) = \{(a_t)_{\ell}\}$. In particular, for any $\bar{a}\in\mathcal{A}_{\ell}$ such that $\bar{a}\neq (a_t)_{\ell}$, there exists $\iota>0$ such that
\begin{equation*}
    \sum_{s\in\mathcal{S}}u_{t_{\ell}}(s,\bar{a})\nu_a(s) \leq \sum_{s\in\mathcal{S}} u_{t_{\ell}}(s,a_{t,\ell})\nu_a(s) - \iota,
\end{equation*}
where we denoted $a_{t,\ell} =(a_t)_{\ell}$, which yields, for any $\theta\in\mathsf{P}(\mathcal{A}_{\ell})$,
\begin{equation*}
    \sum_{(s,\bar a)\in\mathcal{S}\times\mathcal{A}_{\ell}}u_{t_{\ell}}(s,\bar a)\nu_a(s)\theta(\bar a)\leq \sum_{s\in\mathcal{S}} u_{t_{\ell}}(s,a_{t,\ell})\nu_a(s) - \iota (1-\theta(a_{t,\ell})).
\end{equation*}
Moreover, since $a_t\in\mathcal{A}_t(\nu_a^{\star})$, we also have
\begin{equation*}
    \sum_{(s,\bar{a})\in\mathcal{S}\times\mathcal{A}_{\ell}}u_{t_{\ell}}(s,\bar{a})\nu_a^{\star}(s)\theta(\bar{a})\leq \sum_{s\in\mathcal{S}} u_{t_{\ell}}(s,a_{t,\ell})\nu_a^{\star}(s).
\end{equation*}
Combining these two inequalities, we obtain, for any $\lambda_{\ell}\in\mathsf{P}(\mathcal{A}_{\ell})$,
\begin{align*}
    &\sum_{(s,\bar a)\in\mathcal{S}\times\mathcal{A}_{\ell}}u_{t_{\ell}}(s,\bar a)\nu_a^{\epsilon}(s)\theta(\bar a) - \varepsilon\vartheta(\theta,\lambda_{\ell}) 
    \\ &\hspace{5em}\leq \sum_{s\in\mathcal{S}}u_{t_{\ell}}(s,a_{t,\ell})\nu_a^{\varepsilon}(s) - \alpha \iota(1-\theta(a_{t,\ell})) - \varepsilon\vartheta(\theta,\lambda_{\ell}).
\end{align*}
It follows that
\begin{align*}
    &\sum_{s\in\mathcal{S}} u_{t_{\ell}}(s,a_{t,\ell})\nu_a^{\varepsilon}(s) - \varepsilon\vartheta(\delta_{a_{t,\ell}},\lambda_{\ell}) 
    \\ &\hspace{5em}\leq \sum_{s\in\mathcal{S}} u_{t_{\ell}}(s,a_{t,\ell})\nu_a^{\varepsilon}(s) - \alpha \iota(1-\theta^{\star,\varepsilon}_{t_{\ell},\nu_a^{\varepsilon}}(a_{t,\ell})) - \varepsilon\vartheta(\theta^{\star,\varepsilon}_{t_{\ell},\nu_a^{\varepsilon}},\lambda_{\ell}),
\end{align*}
which yields
\begin{equation*}
    0\leq -\iota(1-\theta^{\star,\varepsilon}_{t_{\ell},\nu_a^{\varepsilon}}(a_{t,\ell})) - \varepsilon^{1/2}(\vartheta(\theta^{\star,\varepsilon}_{t_{\ell},\nu_a^{\varepsilon}},\lambda_{\ell}) - \vartheta(\delta_{a_{t,\ell}},\lambda_{\ell})).
\end{equation*}
We know that, up to a subsequence, $(\theta^{\star,\varepsilon}_{t_{\ell},\nu_a^{\varepsilon}})_{\varepsilon>0}$ converges in $\mathsf{P}(\mathcal{A}_{\ell})$ as $\varepsilon\to0$. Assume that the limit is not $\delta_{a_{t,\ell}}$. Then, letting $\varepsilon\to0$ in the previous inequality leads to $0\leq -\iota$, which is impossible. We conclude that
\begin{equation*}
    \theta^{\star,\varepsilon}_{t_{\ell},\nu_a^{\varepsilon}} \underset{\varepsilon\to0}{\longrightarrow} \delta_{a_{t,\ell}}\quad\mbox{and}\quad \theta^{\star}_{\nu_a^{\varepsilon}} \underset{\varepsilon\to0}{\longrightarrow}\delta_{a}.
\end{equation*}
In the end, we obtain that
\begin{align*}
    \int_{\mathsf{P}(\mathcal{S})}W^{\varepsilon}(\nu)d\tau^{\varepsilon}(\nu)  &= (1-\beta)\sum_{a\in\mathcal{R}}p^{\star}_a\sum_{(s,\bar{a})\in\mathcal{S}\times\mathcal{A}^{\mathcal{T}}}w(s,\bar{a})\theta^{\star,\varepsilon}_{\nu_a^{\varepsilon}}(\bar{a})\nu_a^{\varepsilon}(s) 
    \\&\hspace{1em}+ \beta \sum_{(s,\bar{a})\in\mathcal{S}\times\mathcal{A}^{\mathcal{T}}}w(s,\bar{a})\theta^{\star,\varepsilon}_{\bar{\nu}}(\bar{a})\bar{\nu}(s)
    \\ &\underset{\varepsilon\to0}{\longrightarrow}\sum_{a\in\mathcal{R}}p^{\star}_a w(s,a)\nu_a^{\star}(s) = \int_{\mathsf{P}(\mathcal{S})}W(\nu)d\tau^{\star}(\nu),
\end{align*}
since $\nu^{\varepsilon}_a\to \nu_a$ and $\beta\to0$ as $\varepsilon\to0$.

\subsection{Proof of Corollary \ref{cor:convergWeg}}

By Lemma \ref{lem:convergtaureg}, we know that there exists a sequence $(\tau^{\varepsilon})_{\varepsilon>0}\subset\mathsf{T}_{\mu}$ such that
\begin{equation*}
    \bar{W}^{\star} = \underset{\varepsilon\to 0}{\mathrm{liminf}}\;\int_{\mathsf{P}(\mathcal{S})} W^{\varepsilon}(\nu)d\tau^{\varepsilon}(\nu) \leq \underset{\varepsilon\to 0}{\mathrm{liminf}}\;\bar{W}^{\star,\varepsilon} \leq \underset{\varepsilon\to 0}{\mathrm{limsup}}\;\bar{W}^{\star,\varepsilon}.
\end{equation*}
Thus, by Corollary \ref{cor:suboptregtau}, this yields
\begin{equation*}
    \bar{W}^{\star} = \lim_{\varepsilon\to0}\bar{W}^{\star,\varepsilon}.
\end{equation*}

\subsection{Proof of Theorem \ref{theo:conv:regul}}

Let $(\tau^{\varepsilon})_{\varepsilon>0}$ be a sequence of maximizers of \eqref{eq:senderprobreg:tau} belonging in $\mathsf{T}_{M,\mu}$, for some $M\in\mathbb{N}$.
By Lemma \ref{lem:compactau}, we know that this sequence converges, up to a subsequence, to an element  $\tau\in\mathsf{T}_{M,\mu}$. Furthermore, by Corollary \ref{cor:convergWeg} and Lemma \ref{lem:limsupvaluefun}, we have
\begin{equation*}
    \bar{W}^{\star} = \lim_{\varepsilon\to0}\bar{W}^{\star,\varepsilon} = \underset{\varepsilon\to0}{\mathrm{limsup}}\int_{\mathsf{P}(\mathcal{S})}W^{\varepsilon}(\nu)d\tau^{\varepsilon}(\nu)\leq \int_{\mathsf{P}(\mathcal{S})}W(\nu)d\tau(\nu) \leq \bar{W}^{\star}.
\end{equation*}
It follows that $\bar{W}^{\star} = \int_{\mathsf{P}(\mathcal{S})}W(\nu)d\tau(\nu)$ and, thus, $\tau$ is a maximizer of \eqref{eq:senderprob:tau}.

\bibliographystyle{plain}
\bibliography{bibliography.bib}{}

\end{document}